\documentclass[12pt]{article}
\usepackage{amssymb}
\usepackage{latexsym}
\usepackage{amsmath}
\newtheorem{theorem}{Theorem}

\newtheorem{lemma}{Lemma}

\newtheorem{proposition}{Proposition}

\newenvironment{proof}[1][Proof]{\textbf{#1.} }{\ \rule{0.5em}{0.5em}}

\renewcommand{\a}{\mathbf a}
\renewcommand{\r}{\mathbb R}
\renewcommand{\l}{\mathbb L}
\newcommand{\e}{\varepsilon }
\newcommand{\vv}{\varphi }
\newcommand{\A}{\mathcal{A}_{ \alpha, r}(L)}

\def \egal {\stackrel{\mbox{\footnotesize{def}}}{=}}

\begin{document}
\title{Sharp optimality for density deconvolution with dominating bias}
\author{Cristina BUTUCEA$^{1,2}$ and Alexandre B. TSYBAKOV$^1$\\
$^1$Universit\'{e} Paris VI, $^2$Universit\'e Paris X}
\date{August 28, 2004}
\maketitle

\begin{abstract}
We consider estimation of the common probability density $f$ of i.i.d. random
variables $X_i$ that are observed with an additive i.i.d. noise. We assume that the
unknown density $f$ belongs to
a class ${\cal A}$ of densities whose characteristic function is described
by the exponent $\exp(-\alpha |u|^r)$ as $|u|\to\infty$, where $\alpha>0$, $r>0$.
The noise density is supposed
to be known and such that its characteristic function decays as $\exp (-\beta|u|^s)$,
as $|u|\to\infty$, where $\beta>0$, $s>0$. Assuming that $r<s$,
we suggest a kernel type estimator that is optimal in sharp asymptotical minimax sense on ${\cal A}$
simultaneously under the pointwise and the $\l_2$-risks.
The variance of this estimator turns out to be
asymptotically negligible w.r.t. its squared bias.
For $r < s/2$ we construct a sharp adaptive estimator of $f$.
We discuss some effects of dominating bias, such as superefficiency
of minimax estimators.
\end{abstract}



\noindent
{\bf Short title:} Sharp optimality in density deconvolution

%%%%%%%%%%%%%%%%%%%%%%%
\section{Introduction}
%%%%%%%%%%%%%%%%%%%%%%%

Assume that one observes $Y_{1},\ldots ,Y_{n}$ in the
model
\begin{equation*}
Y_{i}=X_{i}+\varepsilon _{i},\quad i=1,\ldots ,n,
\end{equation*}
where $X_{i}$ are i.i.d. random variables with an unknown probability density $f$ w.r.t.
the Lebesgue measure on $\r$,
the random variables $\varepsilon _{i}$ are i.i.d. with known probability density $f^\e$
w.r.t. the Lebesgue mesure on $\r$,
and $(\e_1,\dots, \e_n)$ is independent of
$(X_{1},\ldots ,X_{n})$. The deconvolution problem that we consider here is to estimate $f$
from observations $Y_{1},\ldots ,Y_{n}$.

Denote by $f^{Y} = f*f^\e$ the density of the variables $Y_i$, where $*$ is the convolution sign.
Let $\Phi^{Y}$, $\Phi^{X}$ and $\Phi^{\varepsilon}$
be the characteristic functions of random variables
$Y_{i}$, $X_{i}$ and $\varepsilon_{i}$,
respectively. For an integrable function $g: \r \rightarrow \r$, define
the Fourier
transform
$$\Phi^g(u)=\int{g(x)\exp(ixu)du}.$$

We assume that the unknown density $f$ belongs to the class of functions
\begin{equation*}
\mathcal{A}_{\alpha, r}(L)
= \{ f \ \mbox{is\ a\ probability\ density\ on} \ \r \ \text{and} \, \int{\ \left|
\Phi^f(u) \right|^2 \exp ( 2 \alpha |u|^r )du }\leq 2 \pi L \},
\end{equation*}
where $\alpha >0$, $r>0$, $L>0$ are finite constants.
The  classes of  densities of  this  type have  been studied  by many  authors
starting from  Ibragimov and Hasminskii~(1983).  For a recent overview
see Belitser and Levit~(2001) and Artiles~(2001).

We suppose also in most of the results that the characteristic function of
noise $\e_i$ satisfies the following assumption.

\medskip

{\bf Assumption (N)}.
{\it There exist constants $u_0>0$, $\beta>0$, $s>0$, $b_{\min}>0$, $b_{\max}>0$
and $\gamma, \gamma' \in\r$ such that}
\begin{equation}
\label{N}
b_{\min} |u|^{\gamma}\exp(-\beta \left| u \right|^{s}) \leq \left| \Phi^{\varepsilon}(u)
\right| \leq b_{\max} |u|^{\gamma'}\exp (-\beta \left| u \right|^{s})
\end{equation}
{\it for} $|u| \ge u_0$.

\medskip

Many important probability
densities belong to the class
$\mathcal{A}_{\alpha, r}(L)$ with some $\alpha, r, L$ or
have the characteristic function satisfying (\ref{N}). All such
densities are infinitely many times differentiable on $\r$. Examples include normal,
Cauchy and general stable laws,
Student, logistic, extreme value distributions and other, as well as their mixtures
and convolutions. Note that in these examples the values $r$ and/or $s$ are less or
equal to $2$. Although the densities with $r>2, s>2$ are in principle conceivable, they are
difficult to express in a closed form, and the set of such densities does not contain
statistically famous representatives. This remark concerns especially the noise density
$f^\e$ that should be explicitly known. Therefore, without a meaningful loss, we will
sometimes restrict our study to the case $0<s\le 2$.

For any estimator $\widehat{f}_{n}$ of $f$ define the maximal pointwise
risk over the class $\mathcal{A}_{ \alpha, r}(L)$ for any fixed $x\in \r$ by
\begin{eqnarray*}
R_n(x, \widehat{f}_n, \A)=\sup_{f\in
\mathcal{A}_{ \alpha, r}(L) }E_{f}\left[ \left| \widehat{f}_{n}{(x)-f(x)%
} \right|^{2}\right]
\end{eqnarray*}
and the maximal $\l_2$-risk
\begin{eqnarray*}
R_n(\l_2, \widehat{f}_n, \A)=\sup_{f\in%
\mathcal{A}_{ \alpha, r}(L) }E_{f}\left[ \| \widehat{f}_{n}-f \| _2 ^2%
\right],
\end{eqnarray*}
where $E_{f}(\cdot)$ is the expectation with respect to the joint distribution
$P_{f}$ of $Y_{1},\ldots,Y_{n}$, when the underlying probability density of
$X_{i}$'s is $f$, and $\|\cdot\|_2$ stands for the $\l_2(\r)$-norm. (In what follows
we use the notation $\l_p(\r)$, in general, for the
$\l_p$-spaces of complex valued functions on $\r$.)

\medskip

The asymptotics of optimal estimators differ significantly for the cases $r<s$, $r=s$ and
$r>s$. If $r<s$ the variance of the optimal estimator is
asymptotically negligible w.r.t. the bias, while for $r>s$ the
bias is asymptotically negligible w.r.t. the variance. In this paper we consider the
{\it bias dominated case}, i.e. we assume that $r<s$. The setting with dominating
variance will be treated in another paper.

The problems of density deconvolution with dominating bias were historically the first
ones studied in the literature [cf. Ritov (1987), Stefanski and Carroll (1990), Carroll and Hall (1988),
Zhang (1990), Fan (1991a,b), Masry (1991), Efromovich (1997)],
motivated by the importance of deconvolution with gaussian noise.
These
papers consider, in particular, the noise distributions satisfying (\ref{N}), but the
densities $f$
belonging to finite smoothness classes, such as H\"older or Sobolev ones, where
the estimation of $f$ is harder than for
the class $\A$. In this framework they show that
optimal rates of convergence
are as a power of $\log n$ which suggests that essentially there is no hope to
recover $f$ with a reasonably small error for reasonable sample sizes. This conclusion
is often interpreted
as a general pessimistic message about the gaussian deconvolution problem.
Note, however, that such minimax results are obtained for
the least favorable densities in H\"older or Sobolev classes. Often the underlying density is
much nicer (for instance, it belongs to $\A$,
as the popular densities mentioned above), and the
estimation can be significantly improved, as we show below:
the optimal rates of convergence are in fact faster than
any power of $\log n$.

Pensky~and~Vidakovic~(1999) were  the first  to point out  the effect  of fast
rates in density deconvolution, considering the classes of
densities that are somewhat smaller than $\A$ (including an additional restriction
on the tails of $f$) and with the noise satisfying (\ref{N}). They
analyzed the rates of convergence of
wavelet deconvolution estimators, restricting their attention to the $\l_2$-risk.  Our
results imply that the rates achieved
by their estimators are suboptimal on $\A$ and that the optimal rates can be attained by a
simpler and more traditional kernel deconvolution
method with suitably chosen parameters. We will show that our
method attains not only the optimal rates but also the
best asymptotic constants (i.e.~is sharp optimal).
Moreover, we will prove that the proposed estimator is
sharp optimal simultaneously under the $\l_2$-risk and under the pointwise risk and that it
is sharp adaptive to the parameters $\alpha,r,L$ in some cases.

The most difficult part of our results is the construction of minimax
lower bounds. The technique that we develop
might be useful to get lower bounds for similar ``2 exponents" type settings
in other inverse problems. To our knowledge, except for the case $r=s=1$ treated by
Golubev and Khasminskii (2001), Tsybakov (2000) and Cavalier, Golubev, Lepski and Tsybakov
(2003), such lower bounds are not available even for the Gaussian white noise
(or sequence space) deconvolution model,
although some upper bounds are known (cf. Ermakov (1989),
Efromovich and Koltchinskii (2001)).

Finally, we mention publications on adaptive deconvolution under Assumption (N) or its analogs.
They deal with the problems that
are somewhat different from ours. Efromovich (1997) considered the problem of
deconvolution where the densities $f$ and $f^\varepsilon$ are both periodic
on $[0,2\pi]$, $f^\varepsilon$ satisfies an analog of Assumption (N) expressed in terms
of Fourier coefficients and $f$ belongs to a class of periodic functions of Sobolev type.
He proposed sharp adaptive estimators with logarithmic rates which are optimal for
that framework, as discussed above.
Adaptive deconvolution
in a gaussian white noise model had been studied
by Goldenshluger (1998). He
worked under the Assumption (N) on the Fourier transform of the convolution
kernel or under the assumption that
it decreases as a power of $u$, as $|u|\to \infty$, but he assumed that the function $f$
to estimate belongs to a Sobolev class with unknown parameters.
He proposed a rate adaptive estimator under the pointwise risk.

\section{The estimator, its bias and variance}

Consider  the following kernel estimator of $f$:
\begin{equation}
\hat f_{n}(x)=\frac{1}{nh_{n}}\sum_{i=1}^{n} K_{n}\left( \frac{x-Y_{i}}{h_{n}}
\right), \label{K}
\end{equation}
where $h_n>0$ is a bandwidth and $K_n$ is the function on $\r$ defined as the inverse
Fourier transform of
\begin{equation}
\Phi^{K_{n}}(u)=\frac{I(|u| \le 1)}{\Phi^{\varepsilon}(u/h_{n})}.
\label{K2}
\end{equation}
Here and later $I(\cdot)$ denotes the indicator function.
The function $K_n$ is called kernel, but unlike the usual Parzen-Rosenblatt kernels,
it depends on $n$.

For the existence of $K_n$ it is enough that $\Phi^{K_{n}}\in \l_2(\r)$
(and thus $\Phi^{K_{n}}\in \l_1(\r)$). This holds under mild assumptions.
For example, in view of the continuity property
of characterictic functions, the assumption that
$\Phi^\e (u)\neq 0$ for all $u\in\r$ is sufficient to have $\Phi^{K_{n}}\in \l_2(\r)$. Moreover,
the condition $\Phi^{K_{n}}\in \l_2(\r)$ implies that the kernel $K_n$ is real-valued.
In fact, under this condition we have
$\Phi^{K_{n}}(u)= \Phi^{\varepsilon}(-u/h_{n})V_n(u)$ for almost all $u\in\r$, where
$V_n(u)=I(|u| \le 1)/ |\Phi^{\varepsilon}(u/h_{n})|^2$ is an even real-valued function belonging to
$\l_1(\r)$ and $\Phi^{\varepsilon}(-u/h_{n})$ (the complex conjugate of
$\Phi^{\varepsilon}(u/h_{n})$) is the Fourier transform of real-valued
function $t\mapsto h_n f^\varepsilon (-h_n t )$. This implies that $K_n$ is a convolution of two real-valued functions.

The estimator (\ref{K}) belongs to the family of kernel deconvolution
estimators studied in many papers starting from
Stefanski and Carroll (1990), Carroll and Hall (1988) and
Zhang (1990). It can be also deduced from a unified approach to construction of
estimators in statistical inverse problems (Ruymgaart (1993)).

The following proposition establishes upper bounds on the pointwise and the $\l_{2}$
{\it bias terms}, i.e. on the quantities $| E_{f}\hat f_{n}(x)-f(x)|^{2}$
and $\Vert E_{f}\hat f_{n}-f\Vert _{2}^{2}$.
\begin{proposition}
\label{p1}
Let $f\in \A$, $\alpha>0, r>0, L>0$ and assume that $\Phi^{K_{n}}\in \l_2(\r)$ for any $h_n>0$.
Then the squared bias of $\hat f_{n}(x)$ is bounded as follows
\begin{equation*}
\sup_{x\in\r}\left| E_{f}\hat f_{n}(x)-f(x)\right|^{2}\leq
\frac{L}{2\pi \alpha r}
h_{n}^{r-1}\exp \left( -\frac{2\alpha }{h_{n}^{r}}\right)(1+o(1)),
\end{equation*}
as $h_n \to 0$, while the bias term of the $\l_{2}$-risk satisfies
\begin{equation*}
\Vert E_{f}\hat f_{n}-f\Vert _{2}^{2}\leq L\exp \left( -\frac{2\alpha }{h_{n}^{r}}
\right)
\end{equation*}
for every $h_n>0$.
\end{proposition}
\proof For the pointwise bias we have
\begin{eqnarray*}
\left| E_{f}\hat f_{n}(x)-f(x)\right|^{2}
&=&
\left| \left( \frac{1}{h_{n}}
K_{n}\left( \frac{\cdot }{h_{n}}\right) \ast f^{Y}(\cdot )\right) (x)-f(x)\right|
^{2} \\
&=&\left| \frac{1}{2\pi }\int {\left[ \Phi ^{K_{n}}(uh_{n})\Phi ^{Y}(u)-\Phi
^{X}(u)\right] \exp (-iux)du}\right|^{2} \\
&\leq &
\frac{1}{(2\pi )^{2}}\left( \int { I(|uh_{n}|>1)|\Phi^{X}(u)| du}\right)^{2}.
\end{eqnarray*}
Applying the Cauchy-Schwarz inequality and the assumption that
$f$ belongs to $\mathcal{A}_{\alpha ,r}(L)$ we get
\begin{eqnarray}
&& \left| E_{f}\hat f_{n}(x)-f(x)\right|^{2}\nonumber \\
&\leq &
\frac{1}{(2\pi)^2 }
\int_{|u|>1/h_{n}} \exp (-2\alpha |u|^{r})du
\int_{|u|>1/h_{n}}{|\Phi^{X}(u)|^2\exp (2\alpha |u|^{r})du}
\label{seff3}\\
&\leq &\frac{L}{2\pi }\int_{|u|> 1/h_{n}}{\exp (-2\alpha |u|^{r})du}
\nonumber
\end{eqnarray}
which together with Lemma \ref{lemma2} yields the first inequality of the Proposition.
To prove the second inequality, we apply the Plancherel formula and get
\begin{eqnarray}
\Vert E_{f}\hat f_{n}-f\Vert _{2}^{2}
&= &
\left\Vert \frac{1}{h_{n}}
E_{f}K_{n}\left( \frac{\cdot-Y_1 }{h_{n}}\right) -f(\cdot )\right\Vert _{2}^{2} \nonumber \\
&= &
\frac{1}{2\pi }\int {\left| \Phi^{K_{n}}(uh_{n})\Phi^{Y}(u)-\Phi ^{X}(u)
\right| {^{2}}du} \nonumber\\
& = &
\frac{1}{2\pi}\int { I(|uh_{n}|>1)|\Phi^{X}(u)|^2 du}
\nonumber\\
&\leq &\frac{\exp (-2\alpha /h_{n}^{r})}{2\pi }\int_{\left| u\right| >
1/h_{n}}{\left| \Phi ^{X}(u)\right| {^{2}}\exp (2\alpha \left| u\right|
^{r})du}.
\label{seff2}
\end{eqnarray}
\hfill $\Box $

The next proposition gives upper bounds on the pointwise and the $\l_2$
{\it variance terms} defined as
$$ Var_f \hat f_{n}(x)=E_{f}\left[ |\hat f_{n}(x)-E_{f}\hat f_{n}(x)|^{2}\right] \quad
\text{and} \quad
Var_{f,2}\hat f_{n}=E_{f}\left[\Vert \hat f_{n}-E_{f}\hat f_{n}\Vert ^{2}_2\right]$$
respectively.

\begin{proposition}
\label{p2}
Let the left inequality in (\ref{N}) hold
and $\Phi^\e (u)\neq 0, \forall \ u\in\r$.
Then, for any density $f$ such that $\sup_{x\in\r} f(x) \le f^* <\infty$, the pointwise variance of the estimator $\hat f_{n}(x)$
is bounded as follows
\begin{eqnarray}
\sup_{x\in\r}Var_f \hat f_{n}(x)&= &\sup_{x\in\r}E_{f}\left[ |\hat
f_{n}(x)-E_{f}\hat f_{n}(x)|^{2}\right] \nonumber
\\
&\leq & \min\Big(f_*, \frac{4}{\beta s}h_n^{s-1}\Big) \frac{
h_{n}^{s+2\gamma-1}}{2\pi \beta s b_{\min}^2 n}\exp \left(
\frac{2\beta }{h_{n}^{s}} \right)(1+o(1)), \label{pVar}
\end{eqnarray}
as $h_n \to 0$, and, for an arbitrary density $f$, the variance term
of the $\l_{2}$-risk satisfies
\begin{equation}
\label{L2Var}
Var_{f,2}\hat f_{n}=E_{f}\left[\Vert \hat f_{n}-E_{f}\hat f_{n}\Vert ^{2}_2\right]
\leq
\frac{h_{n}^{s+2\gamma-1}}{2\pi \beta s b_{\min}^2 n}
\exp \left( \frac{2\beta }{h_{n}^{s}}\right)(1+o(1))
\end{equation}
as $h_n \to 0$.
\end{proposition}
\proof
%Define
%$$
%K_{2,n}(x) = \frac{1}{h_{n}} K_{n}^2\left(
%\frac{x }{h_{n}}\right).
%$$
For the pointwise variance we obtain two separate bounds and then
take the minimum of them. To get the first bound, we write
\begin{eqnarray}
Var_f \hat f_{n}(x) &= &
\frac{1}{n}E_{f}\left[ \left| \frac{1}{h_{n}}K_{n}\left(
\frac{x - Y_{1}}{h_{n}}\right) -
E_{f} \left[\frac{1}{h_{n}}K_{n}\left(\frac{x - Y_{1}}{h_{n}}\right)\right]
\right| ^{2}\right]
\nonumber
\\
&\leq &
\frac{1}{nh_{n}}\int \frac{1}{h_n} K_{n}^2\left(\frac{x-y}{h_n}\right) f^{Y} (y)dy \nonumber
\\
&\leq &
\frac{f_*}{nh_{n}}\Vert K_{n}\Vert _{2}^{2},
\label{var1}
\end{eqnarray}
where we used the fact that the convolution density $f^Y = f*f^\e$ is
uniformly bounded by $f_*$.
Applying the Plancherel formula and using (\ref{N}) and (\ref{centre}) of Lemma \ref{lemma2}
in the Appendix we get
\begin{eqnarray}
\Vert K_{n}\Vert _{2}^{2} &=&
\frac{h_n}{2\pi }
\int_{|u|\leq 1/h_{n}} |\Phi ^{\varepsilon }(u)|^{-2} du
\nonumber
\\
&\le&
\frac{h_n}{2\pi b_{\min}^2}
\int_{u_0 \le |u| \le 1/h_{n}} |u|^{- 2\gamma} \exp (2\beta |u|^s) du
+
\frac{h_n}{2\pi }
\int_{|u|\leq u_0} |\Phi ^{\varepsilon }(u)|^{-2} du
\nonumber
\\
&\le&
\frac{h_n}{\pi b_{\min}^2}
\int_{0}^{1/h_{n}} u^{- 2\gamma} \exp (2\beta u^s) du + O(h_n)
\nonumber
\\
&=&
\frac{h_{n}^{s+ 2\gamma}}{2\pi b_{\min}^2 \beta s }\exp \left( \frac{2\beta }{h_{n}^{s}}\right)
(1+o(1)), \quad h_n \to 0.
\label{Ksq}
\end{eqnarray}
This and (\ref{var1}) imply the first bound in (\ref{pVar}). For the
second bound we still use the second line in (\ref{var1}) but then
we apply the Plancherel formula in a different way:
\begin{eqnarray*}
Var_f \hat f_n(x)
& \leq & \frac{1}{n} \int \left(\frac{1}{h_n} K_n\left(\frac{x-y}{h_n}
\right)\right)^2 f^Y (y) dy \\
& \leq & \frac{1}{2 \pi n} \int \left|\Phi^{K_{1,n}^2} (u) \bar
\Phi^Y(u) \right| du \leq \frac{1}{2 \pi n} \int
\left|\Phi^{K_{1,n}^2} (u)\right| du,
\end{eqnarray*}
where  $K_{1,n}(\cdot)=K_n(\cdot/h_n)/h_n$ and  $\bar \Phi^Y$ is the
complex conjugate of $\Phi^Y$. Thus, using that
$\Phi^{K_{1,n}}(u)=\Phi^{K_n} (h_n u)$ and then acting similarly to
(\ref{Ksq}) we get
\begin{eqnarray*}
Var_f \hat f_n(x) & \leq & \frac{1}{2 \pi n} \int \left|
\Phi^{K_{1,n}} \ast \Phi^{K_{1,n}}(u)
 \right| du \leq \frac{1}{2 \pi n} \left( \int \left|\Phi^{K_n} (h_n u)\right| du\right)^2\\
& = & \frac{1}{2 \pi n} \left(\int_{|u| \leq 1/h_n}
\left|\Phi^\varepsilon (u)\right|^{-1}
du \right)^2\\
& \le &
\frac{2 h_n^{2 s + 2\gamma -2}}{\pi\beta^2 s^2 b^2_{\min} n} \exp\left(\frac{2
  \beta} {h_n^s}\right)(1+o(1)),
\end{eqnarray*}
which yields the second bound in (\ref{pVar}). Finally,
\begin{eqnarray*}
Var_{f,2}   \hat   f_{n}   &\leq  &\frac{1}{nh_{n}}\int   \int   \frac{1}{h_n}
K_{n}^2\left(\frac{x-y}{h_n}\right) f^{Y} (y)dy dx \\
&= &
\frac{1}{nh_{n}}\Vert K_{n}\Vert _{2}^{2},
\end{eqnarray*}
and in view of (\ref{Ksq}) we obtain (\ref{L2Var}).\hfill $\Box $

Clearly, the bounds of Proposition \ref{p2} can be applied to $f\in\A$ with,
for example,
$$
f_* = \sup_{f\in\A} \sup_{x\in\r} |f(x)|.
$$
This value is finite and can be taken as in Lemma \ref{lemma1} of the
Appendix.

Interestingly, inequality (\ref{pVar}) shows that asymptotics of the
pointwise variance are different for $0<s<1$ and $s>1$, while this
is not the case for the $\l_{2}$ variance term given by
(\ref{L2Var}). Inequality (\ref{pVar}) can be compared to the recent
result of van Es and Uh~(2004). They studied asymptotic pointwise
variance of the same deconvolution kernel estimator in the
particular case of stable noise distributions with $1/3 < s \leq 2$
and also noticed that $s = 1$ marks a change of behaviour. These
effects concerning variance terms will not be crucial in what
follows since we will consider the bias dominated case.

%%%%%%%%%%%%%%%%%%%%%%%%%%%%%%%%%%%%%%
\section{Optimal bandwidths and upper bounds for the risks}
%%%%%%%%%%%%%%%%%%%%%%%%%%%%%%%%%%%%%%

Propositions \ref{p1} and \ref{p2} lead to upper bounds for pointwise and
$\l_2$ risks that can be minimized in $h_n$. In this section we give an asymptotic approximation
for the result of such a minimization assuming that $r<s$.
The corresponding solutions $h_n$ will be called optimal bandwidths. Note that here we consider
only optimization within a given class of estimators, moreover we minimize upper
bounds on the risks and not the exact risks. However, this turns out to be precise enough in
asymptotical sense: in the next section
we will show that the estimator $\hat{f}_n$ with optimal bandwidth
is sharp minimax over all possible estimators.

Decomposition of the mean
squared error of the kernel estimator into bias and variance terms
and application of Propositions~$\ref{p1}$ and~$\ref{p2}$ yields
\begin{eqnarray*}
E_{f}\left[ \left| {\hat f}_{n}{(x)-f(x)}\right|^{2}\right] &=& \left|
E_{f}\hat f_{n}(x)-f(x)\right| ^{2}+Var_{f}\hat f_{n}(x)
\\
&\leq&
\frac{L}{2 \pi \alpha r} h_{n}^{r-1}\exp \left(-\frac{2\alpha}{h_{n}^{r}}
\right)
+ \frac{f_*}{2 \pi \beta s b_{\min}^2} \frac{h_{n}^{s + 2\gamma-1}}{n}
\exp \left(\frac{2\beta}{h_{n}^{s}} \right).
\end{eqnarray*}
We now minimize the last expression in
$h_n$. Clearly, the minimizer $h_n= \tilde h_n$ tends to $0$, as $n \to \infty$.
Taking derivatives with respect to $h_n$ and neglecting the smaller
terms lead us to the equation for optimal bandwidth
\begin{equation}
\frac{L b_{\min}^2}{f_*}n{\tilde h_n}^{-2\gamma}(1+o(1))=\exp \left(\frac{2\alpha}{{\tilde h_n}^r}
+ \frac{2\beta}{{\tilde h_n}^s} \right), \label{eq}
\end{equation}
(asymptotics are taken as ${\tilde h}_n \rightarrow 0$, $n \to \infty$). Taking
logarithms in the above equation we obtain that the optimal bandwidth
$\tilde h_n$ is a solution in $h$ of the equation
\begin{equation}
2\gamma \log h + \frac{2\alpha}{h^{r}}+ \frac{2\beta}{h^{s}}=\log n +C(1+ o(1)), \label{mse}
\end{equation}
Here and in what follows we denote by $C$ constants with values in $\r$ that can be
different on different occasions.
For the bandwidth $h=\tilde h_n$ satisfying $(\ref{eq})$ and $(\ref{mse})$
we can write
\begin{eqnarray*}
\tilde h_n^{r-1}\exp\left(-\frac{2 \alpha}{\tilde h_n^r}\right)
&=& C(1+o(1))\frac{\tilde h_n^{r + 2\gamma-1}}{n}\exp\left(\frac{2 \beta}{\tilde h_n^s}
\right)
\\
&=& C(1+o(1))\tilde h_n^{r-s}\frac{\tilde h_n^{s + 2\gamma-1}}{n}
\exp\left(\frac{2 \beta}{\tilde h_n^s}\right),
\end{eqnarray*}
with some constant $C>0$.
This proves that, for the optimal bandwidth, the
bias term dominates the variance term
whenever $r<s$. (Strictly speaking, here we consider upper bounds on
the bias and variance terms and not precisely these terms.)

Similarly, for the $\l_2$-risk we get
\begin{eqnarray*}
E_{f}\left[ \| {\hat f}_{n} - f \|_2^2\right] &=& \|
E_{f}{\hat f}_{n} - f \|_2^2+Var_{f,2}{\hat f}_{n}\\
& \leq & L \exp \left(-\frac{2\alpha}{h_{n}^{r}} \right)
+ \frac{1}{2 \pi \beta s b_{\min}^2} \frac{h_{n}^{s + 2\gamma-1}}{n}
\exp \left(\frac{2\beta}{h_{n}^{s}} \right),
\end{eqnarray*}
and the minimizer $h_n = h_n(\l_2)$ of the last expression is a solution in $h$
of the equation
\begin{equation}
(r + 2\gamma-1)\log h + \frac{2\alpha}{h^{r}}+ \frac{2\beta}{h^{s}} = \log n + C(1+o(1)).
\label{mise}
\end{equation}
Now, this equation implies
\begin{eqnarray*}
\exp \left(-\frac{2 \alpha}{h_n^r(\l_2)} \right)
& = & C(1+o(1))\frac{h_n^{r + 2\gamma-1}(\l_2)}{n}\exp \left(\frac{2 \beta}{h_n^s(\l_2)} \right)\\
& = & C(1+o(1)) h_n^{r-s}(\l_2) \frac{h_n^{s + 2\gamma-1}(\l_2)}{n}
\exp\left(\frac{2 \beta}{h_n^s(\l_2)}\right),
\end{eqnarray*}
for some constant $C>0$.
This proves that also for the $\l_2$-risk the bias term dominates
the variance term whenever $r<s$.

Thus  we obtain  two  different equations  $(\ref{mse})$  and $(\ref{mise})$  that
define optimal bandwidths for pointwise  and $\l_2$ risks respectively, and in
both cases the bias terms are asymptotically dominating.

In fact, we can obtain the same  results using a single bandwidth defined as
follows. Denote by $h_*=h_*(n)$ the unique solution of the
equation
\begin{equation}
\frac{2 \beta}{h_*^s} + \frac{2 \alpha}{h_*^r} = \log n - (\log \log n)^2, \label{hstar}
\end{equation}
%and consider the badnwidth
%\begin{equation}
%h_* = h_*(n) = (x_*(n,\alpha,r,\beta,s))^{-1/s}
%\end{equation}
(in what follows we will assume w.l.o.g. that $n \geq 3$ to ensure that
$\log n > (\log \log n)^2$).
Lemma~\ref{nine} in the Appendix implies that, both for the pointwise and the
$\l_2$ loss, the bias terms of the estimator $\hat f_n$ with bandwidth $h_*$ given by
$(\ref{hstar})$ are of the same order as those corresponding to bandwidths $\tilde h_n$ and
$h_n(\l_2)$, while the variance terms corresponding to $(\ref{hstar})$ are
asymptotically smaller. Thus, the pointwise risk and the $\l_2$ risk of
the estimator $\hat f_n$ with bandwidth $h_*$ given by $(\ref{hstar})$ are asymptotically of
the same order as those for estimators $\hat f_n$ with optimal bandwidths $\tilde h_n$
and $h_n(\l_2)$ respectively.

Note that, in fact, $h_*$ is better than both bandwidths $\tilde h_n$ and
$h_n(\l_2)$ in the variance terms, but these terms are asymptotically
negligible w.r.t. the bias ones (cf. Lemma~\ref{nine}). Therefore, the improvement does not appear in
the main term of the asymptotics. Note also that the sequence $(\log \log n)^2$
in $(\ref{hstar})$ can be replaced by a sequence satisfying
$b_n=o((\log n)^{1-r/s})$, $b_n/ \log \log n \to \infty$  and
the above argument remains valid (cf. the proof of Lemma~\ref{nine}).

Calculating the upper bounds for bias terms of the estimator $\hat f_n$ with
bandwidth~$(\ref{hstar})$ we  get the following asymptotical  upper bounds for
its pointwise and $\l_2$ risks respectively:
\begin{equation}
\vv_n^2 = \frac{L}{2 \pi \alpha r} h_*^{r-1}
\exp \left(- \frac{2 \alpha}{h_*^{r}}\right)
= \frac{L}{2 \pi \alpha r}
\left(\frac{\log n}{2 \beta}\right)^{(1-r)/s}
\exp \left(- \frac{2 \alpha}{h_*^{r}} \right)
(1+o(1))
\label{ratevar}
\end{equation}
and
\begin{equation}
\vv_n^2(\l_2) = L\exp \left(- \frac{2 \alpha}{h_*^{r}} \right).
\label{rate2var}
\end{equation}

The above remarks can be summarized as follows.
\begin{theorem}
\label{t1} Let $\alpha >0, L>0, 0<r<s<\infty$, let the left inequality
in (\ref{N}) hold and $\Phi^\e (u)\neq 0, \forall \ u\in\r$.
Then the kernel estimator $\hat f_{n}$  with bandwidth
defined by
(\ref{hstar}) satisfies the following pointwise and $\l_2$-risk bounds
\begin{equation}
\label{upper1}
\limsup_{n \to \infty} \sup_{x \in \r} R_n(x, \hat f_n, \A)\vv_n^{-2}
\leq 1,
\end{equation}
\begin{equation}
\label{upper2}
\limsup_{n \to \infty} R_n(\l_2, \hat f_n, \A)\vv_n^{-2}(\l_2)
\leq 1,
\end{equation}
where the rates $\varphi_{n}$ and $\varphi_{n}(\l_2)$ are given in
(\ref{ratevar}) and (\ref{rate2var}).
\end{theorem}

The case $r = 1$ and $s=2$ is of a particular interest. It covers
the  situation where  the noise  density $f^\e$  is gaussian  ($s=2$)  and the
underlying density  $f$ admits the analytic  continuation into a  strip of the
complex  plane  ($r=1$),  as it  is  the  case  for the  statistically  famous
densities mentioned in the introduction.
%The classes with $r<1$ are even larger.
This case is in the zone $r \leq s/2$, where we get the following behaviour
\begin{equation}
\vv_n^2 = \left\{\begin{array}{cc}
\frac{L}{2 \pi \alpha r} \left(\frac{\log n}{2 \beta}\right)^{(1-r)/s}
\exp \left(-2 \alpha \left(\frac{\log n}{2 \beta}\right)^{r/s}\right)(1+o(1)),
& \text{ if } r< s/2,\\
\frac{L}{2 \pi \alpha r} \left(\frac{\log n}{2 \beta}\right)^{(1-r)/s}
\exp \left(-2 \alpha \sqrt{\frac{\log n}{2 \beta}}+\frac{\alpha^2}{\beta}\right)
(1+o(1)),& \text{ if } r=s/2
\end{array}
\right.
\label{phi1}
\end{equation}
and
\begin{equation}
\vv_n^2(\l_2) = \left\{\begin{array}{cc}
L\exp \left(-2 \alpha \left(\frac{\log n}{2 \beta}\right)^{r/s}\right)(1+o(1)),
& \text{ if } r< s/2,\\
L\exp \left(-2 \alpha \sqrt{\frac{\log n}{2 \beta}}+\frac{\alpha^2}{\beta}\right)
(1+o(1)),& \text{ if } r=s/2.
\end{array}
\right.
\label{phi2}
\end{equation}

The bandwidth (\ref{hstar}) depends on the parameters $\alpha,r$ of the
class $\A$ that are not known in practice. However, it is possible to
construct an adaptive estimator that does not depend on these parameters
and that attains the same asymptotic behavior as in Theorem 1 both for
pointwise and $\l_2$ risks when $r< s/2$.
Define the set of parameters
$$
\Theta = \left\{(\alpha,L,r): \alpha >0, L>0, 0<r< s/2 \right\}.
$$
Note that the parameters $s$ and $\beta$ are supposed to be known since
they characterize the known density of noise $f^\varepsilon$.
\begin{theorem}\label{t2} Suppose that
the left inequality
in (\ref{N}) holds and $\Phi^\e (u)\neq 0, \forall \ u\in\r$.
Let $f^\a_n$ be kernel estimator defined in $(\ref{K})$ with
bandwidth $h_n=h_n^\a$ defined by
\begin{equation}
\label{ha}
h_n^\a= \left(\frac{\log n}{2 \beta} -
\sqrt{\frac{ \log n }{2 \beta}} \right)^{-1/s}
\end{equation}
for $n$ large enough so that $\log n/(2 \beta)>1$.
Then, for all $(\alpha,L,r) \in \Theta$,
\begin{equation*}
\limsup_{n \to \infty}
\sup_{x \in \r}
R_n(x, f_n^\a, \A)\vv_n^{-2} \leq 1,
\end{equation*}
and
\begin{equation*}
\limsup_{n \to \infty}
R_n(\l_2, f_n^\a, \A)\vv_n^{-2}(\l_2) \leq 1,
\end{equation*}
where the rates $\vv_n$ and $\vv_n(\l_2)$ are given in $(\ref{ratevar})$ and
$(\ref{rate2var})$ (and, more particularly, satisfy $(\ref{phi1})$ and $(\ref{phi2})$ with
$r<s/2$).
\end{theorem}
\proof
Since $r/s < 1/2$, we have
$- \left(\frac{\log n}{2 \beta} - \sqrt{\frac{ \log n }{2 \beta}} \right)^{r/s}
> - \frac{\beta}{2\alpha}\sqrt{\frac{ \log n }{2 \beta}}$
for $n$ large enough, and thus
$$
\exp \left( -\frac{2\alpha}{(h_n^\a)^r} \right) \ge
\exp \left( -\beta \sqrt{ \frac{\log n}{2\beta } } \right).
$$
On the other hand,
$$
\frac{1}{n} \exp \left(\frac{2\beta}{(h_n^\a)^s} \right) =
\exp \left( -2\beta \sqrt{ \frac{\log n}{2\beta } } \right).
$$
Therefore, the ratio of the bias term of $f_n^\a$ to the variance term of $f_n^\a$
both for the pointwise
risk and for the $\l_2$-risk is bounded from below by
$$
(\log n)^b \exp \left( \beta \sqrt{ \frac{\log n}{2\beta } } \right)
$$
for some $b\in\r$. This expression tends to $\infty$ as $n\to\infty$. Thus, the
variance terms are asymptotically negligible w.r.t. the bias terms. It remains to check that the
bias terms of $f_n^\a$ for both risks are asymptotically bounded by
$\vv_n^2$ and $\vv_n^2(\l_2)$ respectively.

In view of Proposition~\ref{p1}, for $n$ large enough the bias term of $f_n^\a$
for the pointwise risk is bounded from above by
\begin{eqnarray*}
&&\frac{L}{2 \pi \alpha r} (h_n^\a)^{r-1}
\exp \left( -2\alpha \left( \frac{\log n}{2\beta } \right)^{r/s}
 \left[1- \left(\frac{\log n}{ 2\beta}\right)^{-1/2} \right] ^{r/s}
\right)
\\
&&\le
\frac{L}{2 \pi \alpha r} \left( \frac{\log n}{2\beta }\right)^{(1-r)/s}
\exp \left( -2\alpha \left( \frac{\log n}{2\beta } \right)^{r/s}
 + c\left(\frac{\log n}{ 2\beta}\right)^{r/s-1/2}
\right)(1+o(1))
\\
&&= \vv_n^2(1+o(1)),
\end{eqnarray*}
where $c>0$ is a constant and we have used $(\ref{phi1})$ with $r<s/2$ for the
last equality. Similarly, for $n$ large enough the bias term of $f_n^\a$
for the $\l_2$-risk is bounded from above by
\begin{eqnarray*}
&&L \exp \left( -2\alpha \left( \frac{\log n}{2\beta } \right)^{r/s}
 \left[1- \left(\frac{\log n}{ 2\beta}\right)^{-1/2} \right] ^{r/s}
\right)
\\
&&\le
L \exp \left( -2\alpha \left( \frac{\log n}{2\beta } \right)^{r/s}
 + c\left(\frac{\log n}{ 2\beta}\right)^{r/s-1/2}
\right)
= \vv_n^2(\l_2) (1+o(1)),
\end{eqnarray*}
where $c>0$ and we have used $(\ref{phi2})$ with $r<s/2$ for the last equality.
\hfill $\Box$

If  $r=s/2$, adaptation to  $(\alpha, L)$  is still  possible via  a procedure
similar to that of Theorem~\ref{t2}, but it does not attain the
exact constant, as shows the following result. Introduce the set
$$
\Theta_0 = \{ (\alpha, L): 0 < \alpha \leq \alpha_0, L>0\},
$$
where $\alpha_0 > 0$ is a constant.
\begin{theorem}\label{t3}
Suppose that
the left inequality
in (\ref{N}) holds and $\Phi^\e (u)\neq 0, \forall \ u\in\r$.
Let $f_n^\a$ be the kernel estimator defined in $(\ref{K})$ with bandwidth
$h_n = h_n^\a$ defined by
$$
h_n^\a = \left(\frac{\log n}{2 \beta} - \frac{A}{\beta}
\sqrt{\frac{\log n}{2 \beta}}\right)^{-1/s}
$$
where $A> \alpha_0$ and $n$ is large enough so that $\log n/ (2 \beta)>(A/\beta)^2$.
Then for $r=s/2$ and for all $(\alpha,L) \in \Theta_0$,
\begin{eqnarray}
\limsup_{n\to \infty} \sup_{x \in \r} R_n(x,f_n^\a,\A) \vv_n^{-2} &\leq&
\exp\left(\frac{\alpha A}{\beta} - \frac{\alpha^2}{\beta}\right), \label{ad1}\\
\limsup_{n\to \infty} R_n(\l_2,f_n^\a,\A) \vv_n^{-2}(\l_2) &\leq&
\exp\left(\frac{\alpha A}{\beta} - \frac{\alpha^2}{\beta}\right), \label{ad2}
\end{eqnarray}
where the rates $\vv_n$ and $\vv_n(\l_2)$ are given in $(\ref{phi1})$ and
$(\ref{phi2})$.
\end{theorem}
%
%
%\begin{remark}
%Since  $\alpha \le \alpha_0 <A$,
%the procedure of Theorem~\ref{t3} does not attain
%the exact asymptotical constant: it is adaptive only in the rate of convergence.
%\end{remark}

\noindent \proof It is easily checked that the bias exponent
$$
\exp \left(-\frac{2 \alpha}{(h_n^\a)^r} \right)
=\exp \left(-2 \alpha\sqrt{\frac{\log n}{2 \beta}}
+ \frac{ \alpha A}{\beta}\right)(1+o(1)),
$$
while for the variance term exponent
$$
\frac{1}{n}\exp \left(-\frac{2 \beta}{(h_n^\a)^s} \right)
= \exp\left( - 2A \sqrt{\frac{\log n}{2 \beta}}\right).
$$
Since $A>\alpha$, the bias term of $f_n^\a$
asymptotically dominates its variance term.
Inequalities $(\ref{ad1})$ and $(\ref{ad2})$ now follow from these remarks and
the expressions for $\vv_n^2$, $\vv_n^2(\l_2)$ in $(\ref{phi1})$ and $(\ref{phi2})$
with $r=s/2$.
\hfill $\Box$

%%%%%%%%%%%%%%%%%%%%%%%%%%%%%%
\section{Minimax lower bounds, sharp optimality and superefficiency}
%%%%%%%%%%%%%%%%%%%%%%%%%%%%%%

In this section we establish lower bounds for the risks showing that, under mild additional assumptions, the upper bounds of the previous
section cannot be improved (in a minimax sense on the class of densities $\A$) not only among
kernel estimators, but also among all estimators.
In other words, the estimators suggested in the previous section attain optimal
rates of convergence on $\A$ with optimal exact constants.

We suppose that the following assumption holds.

\smallskip

\noindent {\bf Assumption (ND).}
{\it There exist constants $u_1>0$, $B>0$ and $ \gamma_1 \in \r$ such that $\Phi^\e (u)$
is twice continuously differentiable for $|u| \geq u_1$ with the derivatives satisfying
$$
\max \{|(\Phi^\e (u))^\prime|, |(\Phi^\e(u))^{\prime\prime}|\}
\leq B |u|^{\gamma_1} \exp(-\beta |u|^s),
$$
where $\beta >0$ and $s>0$ are the same as in Assumption~(N).
}

\smallskip

Note that this assumption is satisfied for the examples of popular noise densities mentioned in the
Introduction.

\begin{theorem}
\label{lowbo} Let $\alpha>0, L>0, 0<r< s\le 2$, and suppose that Assumption~(ND) and the right
hand inequality in (\ref{N}) hold. Then
\begin{equation}
\label{lower1}
\liminf_{n \to \infty} \inf_{T_n} R_n(x, T_n, \A)\vv_n^{-2}
\ge 1, \quad \ \forall \ x\in\r,
\end{equation}
and
\begin{equation}
\label{lower2}
\liminf_{n \to \infty} \inf_{T_n}R_n(\l_2, T_n, \A)
\vv_n^{-2}(\l_2)
\ge 1,
\end{equation}
where $\inf_{T_n}$ denotes the infimum over all estimators and the rates
$\vv_n$, $\vv_n(\l_2)$ are defined in
(\ref{ratevar}) and (\ref{rate2var}).
\end{theorem}

Proof of Theorem 4 is given in Section 5.

Theorems 1,2 and 4 immediately imply the following result on
sharp asymptotic minimaxity of the estimators constructed in Section 3.

\begin{theorem}
\label{cor1}
Let $\alpha>0, L>0, 0<r< s\le 2$, let Assumptions~(N), (ND) hold
and $\Phi^\e (u)\neq 0, \forall \ u\in\r$.
Then the kernel estimator $\hat f_{n}$  with bandwidth
defined by (\ref{hstar}) (or with bandwidth defined by (\ref{ha}) if $r<s/2$) is
sharp asymptotically minimax on $\A$ both in pointwise and in $\l_2$ sense:
\begin{equation}
\label{exact1}
\lim_{n \to \infty} R_n(x, \hat f_{n}, \A)\vv_n^{-2}
=
\lim_{n \to \infty} \inf_{T_n} R_n(x, T_n, \A)\vv_n^{-2}
= 1, \quad \ \forall \ x\in\r,
\end{equation}
\begin{equation}
\label{exact2}
\lim_{n \to \infty} R_n(\l_2, \hat f_{n}, \A)\vv_n^{-2}(\l_2)
=
\lim_{n \to \infty} \inf_{T_n}R_n(\l_2, T_n, \A)
\vv_n^{-2}(\l_2)
= 1.
\end{equation}
\end{theorem}

This is the main result of the paper. It shows that the kernel estimator $\hat f_{n}$
with a properly chosen bandwidth $h_n$ is sharp optimal in asymptotically minimax sense
on $\A$ and that for $r<s/2$ the estimator $f_n^\a$ is sharp adaptive in asymptotically
minimax sense on $\A$. Sharp adaptation is thus obtained by direct tuning of the
smoothing parameter without any additional adaptation rule. This is one of the effects
of dominating bias. Theorem \ref{cor1} also
provides exact asymptotical expressions for minimax risks on $\A$ under the
pointwise and the $\l_2$ losses: it states that they are equal to $\vv_n^{2}$ and $\vv_n^{2}(\l_2)$
respectively.

Thus, $\vv_n^{2}$ and $\vv_n^{2}(\l_2)$ can be chosen as reference
values to determine efficiency of estimators. An interesting question is whether
there exist superefficient estimators $\tilde f_{n}$, i.e. such that
\begin{equation}
\label{seff}
\sup_{x \in \r}E_{f}\left[ | \tilde{f}_{n}(x)-f(x) |^2\right] = o(\vv_n^{2}) \quad
\text{ and } \quad
E_{f}\left[ \| \tilde{f}_{n}-f \| _2 ^2\right] = o(\vv_n^{2}(\l_2)),
\end{equation}
as $n\to\infty$,
for any fixed $f\in\A$. The answer to this question is positive, as shows the next proposition.

\begin{proposition}
\label{p3}
Let the conditions of Theorem \ref{t1} hold. Let $\tilde f_{n}$
be the kernel estimator $\hat f_{n}$  with bandwidth
defined by (\ref{hstar}) (or by (\ref{ha}) if $r<s/2$). Then
$\tilde f_{n}$ satisfies
(\ref{seff}). If, moreover, the conditions of Theorem \ref{cor1} hold,
$\tilde f_{n}$ is superefficient in the sense that
\begin{eqnarray}
\label{seff4}
&&\lim_{n\to\infty} \frac{E_{f}[ | \tilde{f}_{n}(x)-f(x) | ^2]}
{\inf_{T_n} \sup_{f\in\A} E_{f}\left[ | T_{n}(x)-f(x) | ^2\right]} = 0, \quad
\forall x \in \r, \\
\label{seff1}
&&\lim_{n\to\infty} \frac{E_{f}[ \| \tilde{f}_{n}-f \| _2 ^2]}
{\inf_{T_n} \sup_{f\in\A} E_{f}\left[ \| T_{n}-f \| _2 ^2\right]} = 0.
\end{eqnarray}
\end{proposition}

{\bf Proof.} Consider the kernel estimator $\hat f_{n}$  with bandwidth
defined by (\ref{hstar}). Instead of using Proposition~\ref{p1} to bound the bias term,
we apply directly (\ref{seff3}) for the pointwise risk and (\ref{seff2}) for the $\l_2$-risk
which yields that, for any fixed $f\in\A$,
\begin{eqnarray*}
\sup_{x \in \r}| E_{f} \hat{f}_{n}(x)-f(x) | ^2 & = &
o\left(h_*^{r-1} \exp (-2\alpha /h_{*}^{r})\right)= o(\vv_n^{2}), \\
\| E_{f} \hat{f}_{n}-f \| _2 ^2 & = &
o\left(\exp (-2\alpha /h_{*}^{r})\right)= o(\vv_n^{2}(\l_2)),
\end{eqnarray*}
as $n\to\infty$.  Now, Proposition \ref{p2} and (\ref{B})  of Lemma \ref{nine}
imply that the variance terms are also $o(\vv_n^2)$ and $o(\vv_n^{2}(\l_2))$, as
$n\to\infty$, respectively.
Hence, (\ref{seff}) follows and implies (\ref{seff4}) and
(\ref{seff1}), in view of Theorem~\ref{cor1}.
The case where the bandwidth is defined by (\ref{ha})
and $r<s/2$ is treated similarly. $\Box$

The result of Proposition \ref{p3} is explained by the fact that the value of the minimax risk in the denominator
of (\ref{seff1}) is attained (up to a $1+o(1)$ factor) on the densities that depend on $n$,
while in the numerator we have a fixed density $f$. Such a superefficiency property occurs
in other nonparametric problems (see e.g. Brown, Low and Zhao~(1997) or Tsybakov~(2004),
Chapter 3), where it is proved for various adaptive estimators. On the contrary, non-adaptive
asymptotically minimax estimators, for example, the Pinsker estimator which is
efficient for ellipsoids in
gaussian sequence model, are not superefficient and turn out to be inadmissible (Tsybakov~(2004),
Section 3.8). Compared with that, the result of Proposition \ref{p3} is somewhat surprising,
because it states that a non-adaptive asymptotically minimax estimator $\hat f_{n}$ with bandwidth
defined by (\ref{hstar}) is superefficient. This provides a simple counter-example
of a superefficient nonparametric estimator which is not adaptive. We conjecture that this is a general property
of nonparametric problems with dominating bias.

\section{Proof of Theorem \ref{lowbo}}

%The rest of the section is devoted to the proof of Theorem \ref{lowbo}.

%%%%%%%%%%%%%%%%%%%%%%%%%%%%%%%%%%%%%%%%
\subsection{General scheme of the proof} \label{5.1}
%%%%%%%%%%%%%%%%%%%%%%%%%%%%%%%%%%%%%%%%

We use the method of proving lower bounds by reduction
to the problem of testing two simple hypotheses (cf. e.g. Tsybakov~(2004), Chapter~2).
Namely, we define two properly chosen probability densities $f_{n1}$ and $f_{n2}$,
depending on $n$ and belonging to $\A$ and we bound the minimax risk as
follows
\begin{eqnarray}
\inf_{{T}_{n}}R_{n}({T}_{n},\mathcal{A}_{\alpha ,r}) \psi_n^{-2}
&\geq &\inf_{{T}_{n}}\max_{f\in \{f_{n1},f_{n2}\}}
E_{f} d^2({T}_{n},f) \psi_n^{-2}\nonumber \\
&\geq &\inf_{T _{n}}\max_{f \in \{f_{n1},f_{n2}\}}
\left(E_f d(T _{n},f) \right)^2 \psi _{n}^{-2} ,
\label{LB0}
\end{eqnarray}
where $R_n(T_n,\A)$ is either $R_n(x,T_n,\A)$ or $R_n(\l_2,T_n,\A)$, $\psi_n$ is
defined as $\vv_n$ or $\vv_n(\l_2)$ (cf. $(\ref{ratevar})$ and $(\ref{rate2var})$)
respectively and $d(T_n,f)$ stands for the distance $|T_n(x)-f(x)|$ at a fixed
point $x$ or the $\l_2$-distance $\|T_n - f\|_2$ respectively. Hence, to prove
the theorem it remains to show that
\begin{equation}
R \egal \inf_{T_n} \max_{f \in \{f_{n1},f_{n2}\}}E_f d(T_n,f)
\geq \psi_n(1+o(1)),\label{LB1}
\end{equation}
as   $n    \to   \infty$,   for   both   pointwise    and   $\l_2$   distances
$d(\cdot,\cdot)$. This will be done  by application of Lemma~\ref{khi2} of the
Appendix.  According to  Lemma~\ref{khi2}, $(\ref{LB1})$  is satisfied  if the
functions $f_{n1}$ and $f_{n2}$ are chosen such that
\begin{eqnarray}
d(f_{n1},f_{n2})  &  \geq  & 2  \psi_n  (1+o(1)),  \text{  as }n  \to  \infty,
\label{LB2}\\
\chi^2(P_{f_{n1}},P_{f_{n2}}) & = & o(1), \text{ as }n \to \infty, \label{LB3}
\end{eqnarray}
where $\chi^2(P_{f_{n1}},P_{f_{n2}})$  is the $\chi^2$-divergence  between the
probability measures $P_{f_{n1}}$ and  $P_{f_{n2}}$ (recall that $P_f$ denotes
the  joint distribution  of $Y_1,\ldots,Y_n$  when the  underlying probability
density of $X_i$'s is $f$).
Thus,  to prove  Theorem~\ref{lowbo} it  suffices to  construct  two functions
$f_{n1}$ and $f_{n2}$ belonging to $\A$ and satisfying $(\ref{LB2})-(\ref{LB3})$.
Since $P_{f_{nj}}$ is  a product of $n$ identical  probability measures corresponding
to the density $f_{nj}^Y= f_{nj}\ast f^\e$, for $j=1,2$, we have
$\chi^2(P_{f_{n1}},P_{f_{n2}})\leq Cn \chi^2(f_{n1}^Y,f_{n2}^Y)$
if $\chi^2(f_{n1}^Y,f_{n2}^Y) \leq 1/n$, where $C$ is a finite constant and
$$
\chi^2(f_{n1}^Y,f_{n2}^Y) = \int \frac{(f_{n1}^Y - f_{n2}^Y)^2}{f_{n1}^Y}(x)dx
$$
(cf. e.g. Tsybakov~(2004), p.~72). Therefore, $(\ref{LB3})$ follows from
\begin{equation}
n \chi^2(f_{n1}^Y,f_{n2}^Y) \to 0, \text{ as } n \to \infty.
\label{LB4}
\end{equation}
We  now proceed to  the construction  of densities  $f_{n1}$, $f_{n2}  \in \A$
satisfying $(\ref{LB4})$ and  $(\ref{LB2})$ for pointwise and $\l_2$-distances
$d(\cdot,\cdot)$.

Consider a density $f_{0}$ of a symmetric stable law whose
characteristic function is
$$
\Phi _{0}\left( u\right) = \left\{\begin{array}{cc}
\exp \left(-\left| c_0 u\right| ^{r}\right) , & \text{ if } 1<r<2,\\
\exp \left(-\left| c_0 u \right| \right), & \text{ if } 0< r \leq 1,
\end{array}
\right.
$$
where $c_0  > \max\{\alpha ^{1/r}, \alpha\}$. Clearly, for any  $0<a<1$ there exists
$c_0>0$ large enough so that $f_0 \in \mathcal{A}_{\alpha,r}(a^2 L)$.
In view of Lemma~$\ref{lemma3}$, there exists $c_1^\prime >0$ such that
\begin{equation}
f_{0}\left( x\right) =\frac{1}{c_0}p\left(\frac{x}{c_0}\right)
\geq \frac{c_1^\prime }{\left| x\right| ^{\max \{r+1, 2\}}+1},
\label{LB5}
\end{equation}
for all $x \in \r$, where $p$ is
the density of stable symmetric distribution with characteristic
function $\exp (-|t|^{\max \{r, 1\}})$, $0<r<2$.
Let $h_+ = h_+(n)$ be the unique solution of the equation
\begin{equation}
\frac{2 \alpha}{h_+^r} + \frac{2 \beta}{h_+^s} = \log n + (\log \log n)^2.
\label{LB6}
\end{equation}
Note that  $h_+$ is analogous to $h_*$  defined by $(\ref{hstar})$
with the only difference that the $(\log \log n)^2$ term changes the sign.

We define the densities $f_{n1}$ and $f_{n2}$ by their characteristic functions
\begin{equation}
\Phi _{n1}\left( u\right) =\Phi _{0}\left( u\right) + \Phi ^{H}\left(
u,h_+\right) , \quad \Phi _{n2}\left( u\right) =\Phi _{0}\left( u\right) - \Phi ^{H}\left(
u,h_+\right), \quad u\in \r, \label{lb1}
\end{equation}
where $u\mapsto \Phi^H(u,h)$ with $h>0$ will be called {\it perturbation  function} and  will be
defined differently  for the pointwise  distance and the  $\l_2$-distance. The
construction of perturbation functions will be based on the following lemma.

\begin{lemma}
\label{auxG} For any $\delta>0$ and any $D>4 \delta$ there exists a
function $\Phi ^{G}: \r \to [0,1]$ such that
\begin{enumerate}
\item[(i)] $\Phi^G$ is  $3$ times continuously differentiable on  $\r$ and the
  first $3$ derivatives of $\Phi^G$ are uniformly bounded on $\r$,

\item[(ii)] $\Phi^G$ is compactly supported on $\left(\delta ,D - \delta \right) $ and
\begin{equation*}
I\left(2\delta \leq u \leq D-2 \delta \right) \leq \Phi ^{G}\left( u\right)
\leq I\left( \delta \leq u \leq D-\delta \right) ,
\end{equation*}
for all $u \in \r$.
\end{enumerate}
\end{lemma}

%%%%%%%%%%%%%%%%%%%%%%%%%%%%%%%%%%%%

\noindent \proof[Proof of Lemma~\ref{auxG}]
Denote by $J_0$ the 5-fold convolution of the indicator function $I(|u| \leq 1)$
with itself. Let $J: \r \to [0,\infty)$ be a rescaling of $J_0$ such that the support
of $J$ is $(-1,1)$ and $\int J(x) dx =1$.
Then $J_0$ and $J$ are 3 times continuously differentiable on $\r$. For
$\delta >0$ and $D> 4 \delta$ define
$$
\Phi^G(u) = \int _{u-D+3\delta/2}^{u-3\delta/2} \frac{2}{\delta} J\left(
\frac{2 x}{\delta}\right)dx.
$$
Clearly, $\Phi^G$ is 3 times continuously differentiable on $\r$ and
$0 \leq \Phi^G (u) \leq 1$, $\forall u \in \r$. Moreover,
$supp~\Phi^G = (\delta, D - \delta)$ and for any $u \in (2\delta, D-2\delta)$
we have $\Phi^G(u) = \int_{-1}^1 J(x) dx = 1.$
\hfill $\Box$

%%%%%%%%%%%%%%%%%%%%%%%%%%%%%%%%%%%%
\subsection{Lower bound at a fixed point}\label{5.2}
%%%%%%%%%%%%%%%%%%%%%%%%%%%%%%%%%%%%

Without loss of generality, we will prove the lower bound for the distance
$d(f,g)=|f(0)-g(0)|$ at the point $x=0$ (if $x\neq 0$ it suffices to shift
the functions $f_{n1}$ and $f_{n2}$ at $x$).
Define the perturbation function
\begin{equation}
\Phi ^{H}\left( u,h\right) =\sqrt{2\pi \alpha rL}~ h^{\left( 1-r\right)
/2}\exp \left( \frac{\alpha }{h^{r}}\right) \exp \left( -2\alpha \left|
u\right| ^{r}\right) \Phi ^{G}\left( \left| u\right| ^{r}-\frac{1}{h^{r}}%
\right) ,\label{lb2}
\end{equation}
where $\Phi ^{G}$ is a function satisfying the properties given in
Lemma~\ref{auxG} for some $\delta>0$ and $D>4\delta$.

Most of the computations below work when $\Phi^G$ is replaced by an
indicator function of the interval $[0,D]$. However, we obviously need a continuous
perturbation function $\Phi^H$ that satisfies  $\Phi^H(0)=0$ to ensure that $f_{n1}$
and $f_{n2}$ integrate to $1$ and that is smooth enough to allow an appropriate bound
on the $\chi^2$-divergence.

\begin{lemma}
\label{pointw}Let $f_{n1}$ and $f_{n2}$ be the functions defined by their
Fourier transforms $\left( \ref{lb1}\right) $,
$(\ref{lb2})$ with $\Phi^G$ satisfying the properties given in Lemma $\ref{auxG}$. Then
we have the following.
\begin{enumerate}
\item The functions $f_{n1}$ and $f_{n2}$ are probability densities for any $n$ large
enough.

\item The functions $f_{n1}$ and $f_{n2}$ belong to $\mathcal{A}_{\alpha ,r}\left( L\right)$
for $n$ large enough if $c_0>0$ in the definition of $f_0$ large enough.

\item The distance between $f_{n1}$ and $f_{n2}$ at $x=0$ satisfies
$$
\left| f_{n1}\left( 0\right) -f_{n2}\left( 0\right) \right|
\geq 2\varphi_{n}[e^{-4 \alpha \delta} - e^{-2 \alpha (D-2 \delta)}] (1+o(1)),
$$
as $n \to \infty$.

\item The $\chi^2$-divergence $\chi ^{2}\left( f_{n1}^{Y},f_{n2}^{Y}\right)$
satisfies $(\ref{LB4})$.

\end{enumerate}
\end{lemma}
\proof $1.$ Clearly, $\Phi^H(\cdot, h)$ is an even, $3$ times
continuously differentiable function on $\r$ having a compact support.
It is easy to see that the integrals $\int |\Phi^H (u,h)|du$ and
$\int |\partial^3 \Phi^H(u,h)/\partial u^3|du$ are bounded uniformly over
$0<h\le h_0$ for any $h_0>0$.
Integration
by parts yields that the inverse Fourier transform of $\Phi^H(\cdot,h)$ can be
written as
\begin{equation}
H(x,h) \egal \frac{1}{2\pi}\int \cos(xu)\Phi^H (u,h)du
= - \frac{1}{2 \pi x^3}\int \sin(xu) \frac{\partial^3 \Phi^H(u,h)}{\partial u^3}du
\end{equation}
for all $x\in\r$ and $0<h\le h_0$.
Thus, there exists a constant $C_H<\infty$ independent of $n$ and such that
\begin{equation}
|H(x,h_+)| \leq C_H (|x|^{3}+1)^{-1}, \text{ for all } x \in \r.
\label{LB7a}
\end{equation}
Denote by $Dom$ the common support of the functions $\Phi^G(|u|^r - 1/h_+^r)$ and
$\Phi^H(u,h_+)$:
$$
Dom \egal \left\{u: |u|^r - \frac{1}{h_+^r} \in [\delta,D-\delta]\right\}
= \left\{ u: \left(\delta +\frac{1}{h_+^r}\right)^{1/r} \leq |u|
\leq \left(D-\delta +\frac{1}{h_+^r}\right)^{1/r}\right\}.
$$
Using the fact that $\left(\delta +1/h_+^r\right)^{1/r} \to \infty$, as $n \to
\infty$, for any fixed $\delta>0$ and applying $(\ref{fin})$ of
Lemma~\ref{lemma2} in the Appendix, we find
\begin{eqnarray}
\Vert H(\cdot, h_+) \Vert_\infty &\egal& \sup_{x \in \r} |H(x,h_+)|
\leq \frac{1}{2 \pi} \int |\Phi^H(u,h_+)|du \nonumber \\
& \leq & \sqrt{\frac{\alpha r L}{2 \pi}} h_+^{(1-r)/2} \exp \left(\alpha /h_+^r \right)
\int_{Dom} \exp (-2 \alpha |u|^r) du \nonumber \\
& \leq & c h_+^{(r-1)/2} \exp (-\alpha /h_+^r)=o(1), \text{ as } n\to \infty,
\label{LB8}
\end{eqnarray}
where $c>0$ is a finite constant.

Now,  $f_{n1}(x)  = f_0(x)+H(x,h_+)$,  $f_{n2}(x)=f_0(x)  - H(x,h_+)$.  Choose
$A>0$ large enough so that for $|x|>A$ we have $C_H(|x|^3+1)^{-1} < c_1^\prime
(|x|^{\max\{r+1,2\}}+1)^{-1}$ (note  that $\max\{r+1,2\} < 3$).  Then, in view
of $(\ref{LB5})$ and $(\ref{LB7a})$,  $f_{nj}(x)>0$, $j=1,2$, for $|x|>A$. Now,
if  $n$  is   large  enough,  $f_{nj}(x)>0$  also  for   $|x|  \leq  A$  since
$\inf_{|x|\leq A} f_0 (x)>0$ (cf. (\ref{LB5})
%Lemma~\ref{lemma3}
) and $(\ref{LB8})$ holds.

Thus,  $f_{nj}(x)>0$, $j=1,2$, for  all $x  \in \r$  if $n$  is large  enough. It
remains  to note  that  $f_{n1}$ and  $f_{n2}$  integrate to  $1$ since  $\int
H(x,h_+)dx = \Phi^H (0,h_+)=0$ (indeed, $0 \not \in supp ~\Phi^H(\cdot, h_+) =
Dom$).

$2.$ We have, by $\left( \ref{lb2}\right) $ and Lemma $\ref{auxG}$,
\begin{eqnarray*}
&&\int \left| \Phi ^{H}\left( u,h_+\right) \right| ^{2}\exp \left( 2\alpha
\left| u\right| ^{r}\right) du \\
&\leq &2\pi \alpha rLh_+^{1-r}\exp \left( \frac{2\alpha }{h_+^{r}}\right)
\int_{Dom}\exp \left( -2\alpha \left| u\right| ^{r}\right) du \\
&\leq &4\pi \alpha rLh_+^{1-r}\exp \left( \frac{2\alpha }{h_+^{r}}\right)
\int_{(\delta + 1/h_+^{r})^{1/r}}^\infty \exp(-2 \alpha u^r)du.
\end{eqnarray*}
By Lemma~\ref{lemma2},
\begin{eqnarray*}
 \int_{(\delta + 1/h_+^{r})^{1/r}}^\infty \exp(-2 \alpha u^r)du
%&=& \frac{1}{2 \alpha r}\left(\delta +\frac{1}{h_+^r}\right) ^{(1-r)/r}
%\exp\left(-2 \alpha \left(\delta + \frac{1}{h_+^r}\right)\right)(1+o(1))\\
&=& \frac{h_+^{r-1}}{2 \alpha r} \exp\left(-\frac{2 \alpha}{h_+^r}\right)
\exp (-2 \alpha \delta) (1+\delta h_+^r)^{(1-r)/r}(1+o(1)),
\end{eqnarray*}
as $n \to \infty$. We get therefore,
\begin{equation}
\int |\Phi^H(x,h_+)|^2 \exp(2 \alpha |u|^r)du \leq 2 \pi L \exp(-2\alpha \delta)
(1+o(1)),
\label{LB9}
\end{equation}
as $n \to \infty$, for any fixed $\delta>0$. Now, choose
$c_0>0$ in  the definition of  $f_0$ large enough  to guarantee that  $f_0 \in
{\mathcal{A}_{ \alpha, r}}(a^2  L)$ with $a=1-e^{-\alpha  \delta/2}$. This  and $(\ref{LB9})$
imply
\begin{eqnarray*}
\left(\int |\Phi_{nj} (u)|^2 \exp(2 \alpha |u|^r) du\right)^{1/2}
& \leq & \|\Phi_0 (\cdot) \exp(\alpha|\cdot|^r)\|_2
+ \|\Phi^H(\cdot,h_+) \exp(\alpha|\cdot|^r)\|_2\\
& \leq  & (1 -e^{-\alpha  \delta/2}) \sqrt{2 \pi L}+e^{-\alpha
\delta} \sqrt{2
  \pi L}(1+o(1))\\
& \leq & \sqrt {2 \pi L}, \quad j=1,2,
\end{eqnarray*}
for $n$ large enough and any fixed $\delta>0$.

$3.$ Using the left inequality in $(ii)$ of Lemma~\ref{auxG} we get
\begin{eqnarray}
\left| f_{n1}\left( 0\right) -f_{n2}\left( 0\right) \right| ^{2}
&=&\frac{1}{\left( 2\pi \right) ^{2}}\left|
\int \left( \Phi _{n1}\left( u\right) -\Phi_{n2}\left( u\right) \right) du\right| ^{2}
=\frac{4}{\left( 2\pi \right) ^{2}}
\left| \int \Phi ^{H}\left( u,h_+ \right) du\right| ^{2}\nonumber \\
&= &\frac{2\alpha rLh_+^{1-r}}{\pi }\exp \left( \frac{2 \alpha }{h_+^{r}}
\right) \left| \int \exp \left( - 2\alpha |u|^{r}\right)\Phi^G\left( |u|^r -
\frac{1}{h_+^r}\right) du \right|^2 \nonumber \\
&\geq & \frac{2\alpha rLh_+^{1-r}}{\pi }\exp \left( \frac{2 \alpha }{h_+^{r}}
\right) \left| 2 \int_{(2 \delta +1/h_+^r)^{1/r}}^{(D-2\delta +1/h_+^r)^{1/r}}
\exp \left( - 2\alpha u^{r}\right) du \right|^2. \label{LB10}
\end{eqnarray}
By $(\ref{fin})$
of Lemma~\ref{lemma2} in the Appendix,
\begin{eqnarray}
  && \int_{(2 \delta +1/h_+^r)^{1/r}}^{(D-2\delta +1/h_+^r)^{1/r}}
  \exp \left( - 2\alpha u^{r}\right) du \nonumber \\
  &=& \frac{h_+^{r-1}}{2 \alpha r}
  \exp \left( -\frac{2 \alpha }{h_+^{r}} \right)\left[ (1+2 \delta h_+^r)^{(1-r)/r}
  e^{-4 \alpha \delta} (1+o(1))\right. \nonumber \\
  && \left. - (1+ (D-2\delta) h_+^r)^{(1-r)/r} e^{-2\alpha(D-2 \delta)}
  (1+o(1)) \right] \nonumber \\
  &=&\frac{h_+^{r-1}}{2 \alpha r}
\exp \left( -\frac{2 \alpha }{h_+^{r}} \right)[e^{-4 \alpha \delta}
-e^{-2\alpha (D-2\delta)}] (1+o(1)),
\label{LB11}
\end{eqnarray}
as $n \to \infty$. The expression in square brackets here is
positive since $D>4\delta$. Combining $(\ref{LB10})$ and
$(\ref{LB11})$ and using $(\ref{2L.1})$ of Lemma~\ref{two} in the
Appendix together with $(\ref{ratevar})$ we get
\begin{eqnarray*}
  \left| f_{n1}\left( 0\right) -f_{n2}\left( 0\right) \right| ^{2}
  & \geq & 4 \left[\frac{L}{2 \pi \alpha r} h_+^{r-1} \exp\left(-\frac{2 \alpha}{h_+^r}
  \right) \right][e^{-4\alpha \delta} - e ^{-2\alpha (D-2\delta)}]^2(1+o(1)) \\
  &= & 4 \left[\frac{L}{2 \pi \alpha r} h_*^{r-1} \exp\left(-\frac{2 \alpha}{h_*^r}
  \right)\right] [e^{-4\alpha \delta} - e ^{-2\alpha (D-2\delta)}]^2(1+o(1))\\
  &=& 4 \vv_n^2 [e^{-4\alpha \delta} - e ^{-2\alpha (D-2\delta)}]^2(1+o(1)),
\end{eqnarray*}
as $n \to \infty$.

$4.$ Inequalities $(\ref{LB5})$, $(\ref{LB7a})$, $(\ref{LB8})$ and the fact that
$r<2$ imply the existence of a constant $c_2^\prime >0$ independent of $n$ and such that
$$
f_{n1}(x) \geq \frac{c_2^\prime}{|x|^{\max \{r+1,2\}}+1}, \quad \forall x \in \r,
$$
for all $n$ large enough.
Since $f^\varepsilon$ is a probability density, we have
$\int_{-M}^M f^\varepsilon (x)dx \geq 1/2$ for a constant
$M>1$ large enough. Hence,
\begin{eqnarray}
f_{n1}^Y(x) &\geq& \int_{-M}^M f_{n1} (x-y) f^\varepsilon(y)dy
\geq \frac{c_2^\prime}{2} \inf_{|y|\leq M}\left[ \frac{1}{|x-y|^{\max\{r+1,2\}}+1}
\right]\nonumber \\
& \geq & c_3^\prime \min\left\{ \frac{1}{M^{\max\{r+1,2\}}},
\frac{1}{|x|^{\max\{r+1,2\}}} \right\} \label{LB12}
\end{eqnarray}
where $n$ and $M$ are large enough, $c_3^\prime >0$ is independent of $n$, and
the last inequality is obtained by considering separately $|x|\le M$ and $|x|>M$. Thus
\begin{eqnarray}
n \chi^2(f_{n1}^Y,f_{n2}^Y) & = & n\int{\frac{(f_{n2}^Y-f_{n1}^Y)^2(x)}{f_{n1}^Y(x)}dx}
= 4 n \int {\frac{(H*f^\varepsilon)^2(x)}{f_{n1}^Y(x)} dx} \nonumber \\
& \leq & \frac{4}{ c_3^\prime} \left(n M^{\max\{r+1,2\}}
\int_{|x|\leq M} (H \ast f^\e)^2(x)dx \right. \nonumber \\
&& \left. + n
\int_{|x|> M} |x|^{\max\{r+1,2\}}(H \ast f^\e)^2(x)dx \right) \nonumber \\
& \leq& (4M^3/ c_3^\prime) (T_{n1}+ T_{n2}), \label{LB13}
\end{eqnarray}
for $n$ and $M$ large enough, where $H(x)=H(x,h_+)$ for brevity and
\begin{equation}
T_{n1} = n \|H\ast f^\e\|_2^2, \quad T_{n2}=n \int |x|^4 (H\ast f^\e)^2(x)dx.
\label{LB13a}
\end{equation}
Using Plancherel's formula and the right hand inequality in $(\ref{N})$ we get,
for $n$ large enough,
\begin{eqnarray}
\Vert H \ast f^{\varepsilon }\Vert_2 ^{2}
&=&\frac{1}{2\pi }\int {\left| \Phi^{H}(u,h_+) \Phi ^{\varepsilon }(u)\right|^{2}du}
\nonumber \\
&\leq &b_{\max}^2 \alpha r L h_+^{1-r} \exp \left(\frac{2 \alpha}{h_+^r}\right)
\int_{Dom} |u|^{2\gamma'}{\exp (-4\alpha |u|^{r}-2\beta |u|^{s})}{du}
\nonumber \\
& \leq & 2b_{\max}^2 \alpha r L h_+^{1-r} \exp \left(\frac{2 \alpha}{h_+^r}\right)
\int_{(\delta +1/h_+^r)^{1/r}}^{\infty}
u^{2\gamma'}
{\exp (-4\alpha u^{r}-2\beta u^{s})}{du}
\nonumber \\
&\leq & 2 b^2_{\max} \alpha r L h_+^{1-r} \exp \left(-\frac{2 \alpha}{h_+^r} \right)
\int_{1/h_+}^{\infty} u^{2\gamma'} {\exp (-2\beta u^{s})}{du}.
\label{LB14}
\end{eqnarray}
The last integral is evaluated using $(\ref{fin})$ of Lemma~\ref{lemma2} in
the Appendix:
\begin{equation}
\int_{1/h_+}^{\infty} u^{2\gamma'} {\exp (-2\beta u^{s})}{du}
= \frac{h_+^{s-2\gamma'-1}}{2 \beta s} \exp\left(-\frac{2 \beta}{h_+^s} \right) (1+o(1)),
\label{LB14a}
\end{equation}
as $n \to \infty$. This, together with $(\ref{LB14})$ and (\ref{2L.2}) of
Lemma~\ref{two} in the Appendix, yields
\begin{equation}
\|H \ast f^\e \|_2^2 \leq C h_+^{s-2\gamma'-r} \exp\left( -\frac{2 \alpha}{h_+^r}
-\frac{2 \beta}{h_+^s}\right)=o\left( \frac{1}{n}\right),
\label{LB14b}
\end{equation}
as $n \to \infty$, where $C>0$ is a constant.
Thus,
\begin{equation}\label{LB15}
    T_{n1}=o(1), \text{ as } n \to \infty.
\end{equation}

Now, assume that $n$ is large enough to have $(\delta + 1/h_+^r)^{1/r} > \max(u_0,u_1)$,
where $u_0>0$, $u_1>0$ are the constants in Assumptions~(N) and~(ND).
Then $\Phi^G(|u|^r-1/h_+^r)=0$
for $|u| \leq \max(u_0,u_1)$, and thus the function $\Phi^H(\cdot,h_+) \Phi^\e(\cdot)$ is
twice continuously differentiable on $\r$. Using Assumption~(ND), the right hand
inequality in $(\ref{N})$ and the fact that $\Phi^G$, together with its first two
derivatives, is uniformly bounded on $\r$ we find that there exist constants
$B_1 <\infty$ and $a \in \r$ such that, for $n$ large enough and all $u\in\r$,
\begin{equation}\label{LB16}
    \left| (\Phi^H(u,h_+)\Phi^\e(u))^{\prime \prime} \right|
    \leq B_1 h_+^{(1-r)/2} \exp \left(\frac{\alpha}{h_+^r}\right) |u|^a
    \exp(-2\alpha |u|^r - \beta |u|^s).
\end{equation}
Thus, for $n$ large enough, we have, by Plancherel's formula for derivatives and
$(\ref{LB16})$,
\begin{eqnarray}
  T_{n2} &=& \frac{n}{2 \pi}
  \int\left|(\Phi^H(u,h_+)\Phi^\e(u))^{\prime \prime} \right|^2 du\nonumber \\
  &\leq & \frac{n}{2 \pi} B_1^2 h_+^{1-r} \exp\left( \frac{2 \alpha}{h_+^r}\right)
  \int_{Dom} |u|^{2a} \exp(-4\alpha|u|^r - 2\beta |u|^s) du
  \nonumber \\
  &\leq &\frac{n}{\pi} B_1^2 h_+^{1-r} \exp\left( \frac{2 \alpha}{h_+^r}\right)
  \int_{(\delta+1/h_+^r)^{1/r}}^{\infty} u^{2a} \exp(-4\alpha u^r - 2\beta u^s) du
  \nonumber \\
  & \leq & \frac{n}{\pi} B_1^2 h_+^{1-r} \exp\left( -\frac{2 \alpha}{h_+^r}\right)
  \int_{1/h_+}^\infty u^{2a} \exp( - 2 \beta u^s) du.
  \label{LB17}
\end{eqnarray}
%Here, in view of Lemma~\ref{lemma2} in the Appendix,
%\begin{equation}\label{LB17a}
%\int_{1/h_+}^\infty u^{2a} \exp( - 2 \beta u^s) du = \frac{1}{2 \beta s}
%h_+^{-2a +s -1} \exp\left(-\frac{2 \beta}{h_+^s} \right) (1+o(1)),
%\end{equation}
%as $n \to \infty$.
Plugging (\ref{LB14a}) with $\gamma'=a$ into $(\ref{LB17})$ and using (\ref{2L.2}) of
Lemma~\ref{two} in the Appendix we get
\begin{equation}\label{LB18}
T_{n2}\leq C n h_+^{-2a +s -r}
\exp\left(-\frac{2 \alpha}{h_+^r} - \frac{2 \beta}{h_+^s} \right) (1+o(1))
=o(1),
\end{equation}
as $n \to \infty$, where $C>0$ is a constant.

Combining $(\ref{LB13})$, $(\ref{LB15})$ and $(\ref{LB18})$ we get that
$n \chi^2 (f_{n1}^Y, f_{n2}^Y) \to 0$, as $n\to \infty$.
\hfill $\Box$

{\bf Proof of $(\ref{lower1})$.} We use the general scheme of Section~\ref{5.1} with
$d(f_{n1},f_{n2})= |f_{n1}(0) - f_{n2}(0)|$. Choose $c_0>0$
in the definition of $f_0$ large enough to guarantee that assertion 2 of Lemma~\ref{pointw}
holds. Lemma~\ref{pointw}
implies that $(\ref{LB4})$ and thus $(\ref{LB3})$ are satisfied and that
$(\ref{LB2})$ holds with
$$
\psi_n = \vv_n [e^{-4 \alpha \delta} - e^{-2 \alpha (D-2\delta)}].
$$
Therefore, Lemma~\ref{khi2} of the Appendix implies that
$$
R \geq \vv_n[e^{-4 \alpha \delta} - e^{-2 \alpha (D-2\delta)}] (1+o(1)),
$$
as $n \to\infty$, where $R$ is defined in $(\ref{LB1})$. This and $(\ref{LB0})$
yield that, as $n \to \infty$,
$$
\inf_{T_n} R_n(0,T_n,\A) \vv_n^{-2}
\geq [e^{-4 \alpha \delta} - e^{-2 \alpha (D-2\delta)}](1+o(1)).
$$
Taking limits as $n \to \infty$ and then as $D \to \infty$ and $\delta \to 0$
we get $(\ref{lower1})$ for $x=0$. The proof for $x\neq 0$ is analogous (see the remark
at the beginning of this section).
\hfill $\Box$

%%%%%%%%%%%%%%%%%%%%%%%%%%%%%%%%%
\subsection{Lower bound in $\l_2$}
%%%%%%%%%%%%%%%%%%%%%%%%%%%%%%%%%

Introduce the perturbation function
\begin{equation}
\Phi^H(u,h)  =   \sqrt{2  \pi  \alpha  r  L   (d-1)} \ h^{(1-r)/2}  e^{(d-1)\alpha/h^r}
    \exp \left( - \alpha d |u|^r \right)\Phi^G
    \left(|u|^r - \frac{1}{h^r} \right), \label{lb3}
\end{equation}
where $\Phi^G$ is a function satisfying the properties given in Lemma~\ref{auxG}
and $d=d(\delta)>1$ is a constant depending on the value $\delta$ that appears in the
construction of $\Phi^G$.
The argument below is similar to that of Section~\ref{5.2}, modulo the choice of the
perturbation function~$(\ref{lb3})$ which is slightly different from $(\ref{lb2})$.
The argument goes through with $d$ such that $d(\delta) \to \infty$ and
$\delta d(\delta) \to 0$ as $\delta \to 0$, but we will set for simplicity
$d(\delta) = \delta^{-1/2}$ and assume that $0<\delta<1$, which ensures that
$d(\delta)>1$.

\begin{lemma}\label{L2dist}
Let $f_{n1}$ and $f_{n2}$ be the functions defined by their Fourier transforms
$(\ref{lb1})$, $(\ref{lb3})$ with $\Phi^G$ satisfying the properties of
Lemma $\ref{auxG}$ and $0<\delta<1$. Then we have the following.
\begin{enumerate}
\item The functions $f_{n1}$ and $f_{n2}$ are probability densities for $n$
large enough.

\item The functions $f_{n1}$ and $f_{n2}$ belong to $ \A$ for $n$ large
enough if $c_0>0$ in the definition of $f_0$ large enough.

\item The $\l_2$ distance between $f_{n1}$ and $f_{n2}$ satisfies
$$
\| f_{n1} - f_{n2} \|_2 \geq 2\varphi_{n}(\l_2) \left((1- \sqrt{\delta})
[e^{-4 \alpha \sqrt{\delta}} - e^{-2\alpha (D-2 \delta)/ \sqrt{\delta}}])\right)^{1/2}
(1+o(1)),
$$
as $n \to \infty$.

\item  The $\chi^2$-divergence $ \chi ^{2}\left( f_{n1}^{Y},f_{n2}^{Y}\right)$
satisfies $(\ref{LB4})$.
\end{enumerate}
\end{lemma}
\proof $1.$ The argument is analogous to the proof of assertion~1 of Lemma~\ref{pointw}.
In particular, one also has $|H(x,h)| \leq C_H^\prime (|x|^3+1)^{-1}$,
$\forall x \in \r$, and $\|H(\cdot,h_+)\|_\infty =o(1)$, as $n \to \infty$,
for some constant $C_H^\prime < \infty$. We omit the details.

$2.$ We have by $\left( \ref{lb1}\right) $ and Lemma $\ref{auxG}$
\begin{eqnarray*}
&&\int \left| \Phi ^{H}\left( u,h_+\right) \right| ^{2}\exp \left( 2\alpha
\left| u\right| ^{r}\right) du \\
&\leq &2\pi \alpha rL(d-1)h_+^{1-r}\exp \left( \frac{2(d-1)\alpha }{h_+^{r}}\right)
\int_{Dom}\exp \left( -2\alpha (d-1)\left| u\right| ^{r}\right) du \\
&\leq &4\pi \alpha rL(d-1)h_+^{1-r}\exp \left( \frac{2(d-1)\alpha }{h_+^{r}}\right)
\int_{(\delta+1/h_+^r)^{1/r}}^\infty \exp \left( -2\alpha (d-1) u^{r}\right) du.
\end{eqnarray*}
By Lemma~\ref{lemma2},
\begin{eqnarray*}
  && \int_{(\delta+1/h_+^r)^{1/r}}^\infty \exp \left( -2\alpha (d-1) u^{r}\right) du\\
 % &=& \frac{1}{2\alpha (d-1)r} \left(\delta +\frac{1}{h_+^r}\right)^{(1-r)/r}
  %\exp\left( -2\alpha (d-1) \left(\delta +\frac{1}{h_+^r}\right) \right)
  %(1+o(1))\\
  &=& \frac{h_+^{r-1}}{2\alpha (d-1)r} \exp\left(- \frac{2(d-1)\alpha}{h_+^r}\right)
  \exp\left( -2\alpha (d-1) \delta \right)\left(1+\delta h_+^r \right)^{(1-r)/r}
  (1+o(1)),
\end{eqnarray*}
as $n \to \infty$. We get therefore,
$$
\int |\Phi^H (u,h_+)|^2 \exp(2\alpha|u|^r) du
\leq 2 \pi L \exp(-2 \alpha (d-1)\delta) (1+o(1)),
$$
as $n \to \infty$, for any fixed $\delta>0$. Now, since
$d=\delta^{-1/2}$, we get that the last exponent is strictly less than $1$ for
$0<\delta<1$, and thus the argument similar to that after formula $(\ref{LB9})$
can be applied to show that
$$
\int |\Phi_{nj} (u)|^2 \exp(2\alpha|u|^r) du
\leq 2 \pi L, \quad j=1,2,
$$
for $n$ large enough, if $c_0>0$ in the definition of $f_0$ is chosen large enough.

$3.$ The $\l_2$ distance is
\begin{eqnarray}
&&\| f_{n1} -f_{n2} \|_2 ^{2} = \frac{1}{2\pi}
\int \left( \Phi _{n1}\left( u\right) -\Phi_{n2}\left( u\right) \right)^2 du
=\frac{4}{ 2\pi} \int \left| \Phi ^{H}\left( u,h_+\right)\right|^2 du
\nonumber \\
&=& 4 L\alpha r (d-1)h_+^{1-r}\exp \left( \frac{2 (d-1) \alpha }{h_+^{r}}\right)
\int \exp \left( -2 \alpha d \left| u\right| ^{r}\right)
\left| \Phi^G \left(|u|^r-\frac{1}{h_+^r} \right)\right|^2 du\nonumber \\
&\geq &4 L \alpha r (d-1)h_+^{1-r}\exp \left( \frac{2 (d-1) \alpha }{h_+^{r}}
\right) \left[2 \int_{(2\delta + 1/h_+^r)^{1/r}}^{(D-2\delta +1/h_+^r)^{1/r}}
\exp \left( -2 \alpha d u^{r}\right) du\right]
\label{LB19}
\end{eqnarray}
where we used the left inequality in $(ii)$ of Lemma~\ref{pointw}.
Lemma~\ref{lemma2} implies that (cf. $(\ref{LB11})$):
\begin{eqnarray*}
  && \int_{(2\delta + 1/h_+^r)^{1/r}}^{(D-2\delta +1/h_+^r)^{1/r}}
\exp \left( -2 \alpha d u^{r}\right) du \\
  &=& \frac{h_+^{r-1}}{2 \alpha d r} \exp\left(-\frac{2\alpha d}{h_+^r} \right)
  [e^{-4\alpha d \delta} - e^{-2 \alpha d(D-2\delta)}](1+o(1)),
\end{eqnarray*}
as $n \to \infty$. Substituting this into $(\ref{LB19})$ and using
$(\ref{2L.1})$ of Lemma~\ref{two} we obtain
\begin{eqnarray*}
  \| f_{n1} -f_{n2} \|_2 ^{2} &\geq & 4 L\frac{d-1}{d}
  \exp\left(-\frac{2\alpha}{h_+^r} \right)
  [e^{-4\alpha d \delta} - e^{-2 \alpha d(D-2\delta)}](1+o(1)) \\
 % &=& 4 L \exp\left(-\frac{2\alpha }{h_+^r} \right) (1-\sqrt{\delta})
 % [e^{-4\alpha \sqrt{\delta}} - e^{-2 \alpha (D-2\delta)/\sqrt{\delta}}]
 % (1+o(1)) \\
  &=& 4 L \exp\left(-\frac{2\alpha }{h_*^r} \right) (1-\sqrt{\delta})
  [e^{-4\alpha \sqrt{\delta}} - e^{-2 \alpha (D-2\delta)/\sqrt{\delta}}]
  (1+o(1)) \\
  &=& 4 \vv_n^2(\l_2)(1-\sqrt{\delta})[e^{-4\alpha \sqrt{\delta}} - e^{-2 \alpha
  (D-2\delta)/\sqrt{\delta}}]
  (1+o(1)),
\end{eqnarray*}
as $n \to \infty$, (cf. the definition of $\vv_n (\l_2)$ in
$(\ref{rate2var})$).

$4.$ Similarly to the proof of assertion 4 of Lemma~\ref{pointw},
we obtain
\begin{equation}\label{LB20}
    n \chi^2(f_{n1}^Y, f_{n2}^Y) \leq c_4^\prime(T_{n1}+T_{n2}),
\end{equation}
for $n$ and $M$ large enough, where $T_{n1}$ and $T_{n2}$ are
defined in $(\ref{LB13a})$ and $c_4^\prime < \infty$ is a constant. The
only difference from the proof of Lemma~\ref{pointw} is that the
function $H(x)=H(x,h_+)$ is now defined as the inverse Fourier
transform of $(\ref{lb2})$ and not as that of $(\ref{lb1})$. As in
$(\ref{LB14})-(\ref{LB14b})$, we get, for $n$ large enough,
\begin{eqnarray}
&& T_{n1} =  n \|H \ast f^\e\|_2^2 \nonumber \\
&\leq & b_{\max}^2 \alpha r L (d-1)
n h_+^{1-r}\exp\left(\frac{2 (d-1)\alpha}{h_+^r}\right)
\int _{Dom}
|u|^{2\gamma'}
{\exp \left( -2 \alpha d |u|^r -2 \beta |u|^s \right) du}
\nonumber \\
& \leq &c^\prime n h_+^{1-r}\exp\left(-\frac{2\alpha}{h_+^r}\right)
\int _{1/h_+}^\infty
u^{2\gamma'}
{\exp \left( -2 \beta u^s \right) du} \nonumber \\
& \leq &c^{\prime \prime} n h_+^{s-2\gamma'-r}\exp\left(-\frac{2\alpha}{h_+^r}
-\frac{2 \beta}{h_+^{s}}\right)=o(1),
\label{LB21}
\end{eqnarray}
as $n\to \infty$, where $c^\prime>0$ and $c^{\prime \prime}>0$ are some
finite constants.

Next, similarly to $(\ref{LB16})$, we have, for $n$ large enough
and all $u\in\r$,
\begin{eqnarray*}
&& |(\Phi ^H(u,h_+) \Phi^\varepsilon(u))^{\prime \prime}|
\leq  B_2 h_+^{(1-r)/2}\exp\left(\frac{(d-1)\alpha}{h_+^r}\right)
|u|^{a^\prime} \exp(-2\alpha d |u|^r -\beta |u|^s ),
\end{eqnarray*}
where $B_2 <\infty$ and $a^\prime \in \r$ are some constants.
This implies, as in $(\ref{LB17})-(\ref{LB18})$, that
\begin{eqnarray}
  T_{n2}&=&\frac{n}{2 \pi}
  \int |(\Phi ^H(u,h_+) \Phi^\varepsilon(u))^{\prime \prime}|^2 du
  \nonumber \\
  &\leq &\frac{n}{\pi} B_2^2 h_+^{1-r}\exp\left(-\frac{2 \alpha}{h_+^r}\right)
  \int_{1/h_+}^{\infty}u^{2a^\prime} \exp(- 2\beta u^s )du\nonumber \\
  &\leq & \bar{c} n h_+^{-2 a^\prime +s-r} \exp\left(-\frac{2 \alpha}{h_+^r}
  -\frac{2 \beta}{h_+^s}\right)=o(1),\label{LB22}
\end{eqnarray}
as $n \to \infty$, where $\bar{c}>0$ is finite constant. It remains now
to combine $(\ref{LB20})-(\ref{LB22})$.

{\bf Proof of $(\ref{lower2})$} is now obtained following the same lines
as the proof of $(\ref{lower1})$ in Section~\ref{5.2}, but with $d(f_{n1},f_{n2})
=\|f_{n1}-f_{n2}\|_2$ and $\psi_n=\vv_n(\l_2) \Big((1-\sqrt{\delta})
[e^{-4\alpha\sqrt{\delta}}-e^{-2 \alpha(D-2 \delta)/\sqrt{\delta}}]\Big)^{1/2}$.
\hfill $\Box$

%%%%%%%%%%%%%%%%%%%%%%%%%%%%%%%%
\section{Appendix}
%%%%%%%%%%%%%%%%%%%%%%%%%%%%%%%%

Let $({\cal X, A})$ and $(\Theta,{\cal T})$ be measurable spaces and let $P_1$
and $P_2$ be two probability measures on $\cal A$. Let
$d:(\Theta \times \Theta, \cal T \otimes \cal T) \to (\r_{+},{\cal B})$ be a non-negative
measurable function where ${\cal B}$ is the Borel $\sigma$-algebra. Define
\[
R=\inf_{\hat{\theta}} \max_{i \in \{1,2\}} E_i[d(\hat{\theta},\theta_i)],
\]
where $\inf_{\hat{\theta}}$ denotes the infimum with respect to all the measurable mappings
$\hat{\theta}:({\cal X,A}) \to (\Theta, {\cal T})$, $E_i$ denotes the expectation
with respect to $P_i$, and $\theta_1$, $\theta_2$ are two elements of $\Theta$.
\begin{lemma}\label{khi2} Suppose that:
\begin{itemize}
\item[(i)] $d(\cdot,\cdot)$ satisfies the triangle inequality,
\item[(ii)] $\theta_1,\theta_2 \in \Theta$ are such that
$d(\theta_1,\theta_2) \geq 2 \psi,$
for some $\psi >0$,
\item[(iii)] $P_2 \ll P_1$ and there exist constants $\tau>0$ and $0<\gamma_0<1$ such that
\[
P_1 \left[\frac{dP_2}{dP_1} \geq \tau \right]\geq 1-\gamma_0.
\]
\end{itemize}
Then
\begin{equation}
R\geq \psi (1-\gamma_0)\min\{\tau,1\}. \label{exactlb}
\end{equation}
Furthermore, if instead of ${\it (iii)}$ we suppose that
\begin{itemize}
\item[(iv)] $\chi^2 (P_1,P_2) \leq \gamma_0^2$, where $0<\gamma_0<1$ and
\end{itemize}
\[
\chi^2(P_1,P_2)=\int \left( \frac{dP_2}{dP_1}-1 \right)^2 dP_1,
\]
then
\begin{equation}
R\geq \psi (1-\gamma_0)(1-\sqrt{\gamma_0}).\label{exactchi2}
\end{equation}
\end{lemma}
\proof We first show $(\ref{exactlb})$. We have
\begin{eqnarray*}
R&\geq& \frac{1}{2} \inf_{\hat{\theta}} \left( E_1 [d(\hat\theta,\theta_1)]
+ E_2 [d(\hat \theta,\theta_2)]\right)\\
&\geq& \frac{1}{2} \inf_{\hat{\theta}} \left( E_1 [d(\hat\theta,\theta_1)]
+ \tau E_1 \left[I\left(\frac{dP_2}{dP_1}\geq \tau \right)d(\hat \theta,\theta_2)
\right] \right)\\
&\geq& \frac{\min\{\tau,1\}}{2} \inf_{\hat{\theta}}
E_1 \left[I\left(\frac{dP_2}{dP_1}\geq \tau \right)
[d(\hat \theta,\theta_1)+d(\hat \theta,\theta_2) ]\right].
\end{eqnarray*}
Using here the triangle inequality and ${\it (ii)-(iii)}$, we find
\[
R\geq \psi  \min\{\tau,1\} P_1 \left[ \frac{dP_2}{dP_1}\geq \tau\right]
\geq \psi (1-\gamma_0) \min\{\tau,1\}.
\]
To show $(\ref{exactchi2})$ it is sufficent to note that, in view of Chebyshev's
inequality
\[
P_1\left[ \frac{dP_2}{dP_1}\geq 1-\sqrt{\gamma_0}\right] =1-P_1\left[\frac{dP_2}{dP_1}
-1<-\sqrt{\gamma_0} \right] \geq 1-\frac{1}{\gamma_0} \int \left( \frac{dP_2}{dP_1}-1
\right)^2 dP_1 \geq 1-\gamma_0,
\]
and thus ${\it (iv)}$ implies ${\it (iii)}$ with $\tau=1-\sqrt{\gamma_0}$.\hfill $\Box$

%%%%%%%%%

\begin{lemma} \label{lemma1}
For $0<\alpha, r,L<\infty$,
\begin{eqnarray*}
\sup_{f \in {\cal A}_{\alpha,r}(L)} \sup_{x \in \r} |f(x)|
& \leq & L+\pi^{-1} C(r,\alpha),\\
\end{eqnarray*}
where $C(r,\alpha)=\int_0^\infty \exp(-2\alpha u^r) du$.
\end{lemma}
\proof Let $\Phi=\Phi^f $ be the characteristic function of $f$. Clearly,
\begin{equation}\label{A6.1}
|f (x)| \leq \frac{1}{2 \pi} \int |\Phi(u)| du,\quad
\forall x \in \r.
\end{equation}
By Markov's inequality
\begin{eqnarray*}
\int |\Phi(u)| \ I\Big(|\Phi(u)|\exp \left( 2\alpha |u|^r \right)>1\Big)   du & \leq &
\int \exp \left( 2\alpha |u|^r \right) |\Phi(u)|^2 du \leq 2 \pi L.
\end{eqnarray*}
Also,
\begin{eqnarray*}
\int |\Phi(u)| \ I\Big(|\Phi(u)|\exp \left( 2\alpha |u|^r \right)\leq 1\Big)   du & \leq &
2\int_0^\infty \exp \left(- 2\alpha u^r \right) du = 2 C(r,\alpha).
\end{eqnarray*}
Combining the last two inequalities with $(\ref{A6.1})$ proves the Lemma. \hfill $\Box$

%%%%%%%%%%%%%%%%%%%%%%%%%%%%

\begin{lemma} \label{lemma2}
For any positive $\alpha$, $\beta$, $r$, $s$ and for any $A\in\r$ and
$B\in\r$, we have
\begin{equation}
\int_v^{\infty} u^A \exp \left(-\alpha u^r \right)du = \frac{1}{\alpha r}
v^{A+1-r} \exp(-\alpha v^r)(1+o(1)), v \to \infty,
\label{fin}
\end{equation}
and
\begin{equation}
\int_0^v u^B \exp \left(\beta u^s \right)du = \frac{1}{\beta s}
v^{B+1-s} \exp(\beta v^s)(1+o(1)), v \to \infty.
\label{centre}
\end{equation}
\end{lemma}

Proof of this lemma is omitted. It is based on integration by parts and
standard evaluations of integrals.

%%%%%%%%%%%%%%%%%%%%%%%%%%%%%

\begin{lemma} \label{lemma3}
Let $p$ be the density of stable symmetric distribution with characteristic
function $\exp (-|t|^r)$, $1<r<2$. Then $p$ is continuous, $p(x)>0$ for all $x \in \r$ and there
exist $c_1>0$, $c_2>0$ such that
\[
p(x) \geq c_1 |x|^{-r-1},
\]
for $|x|\geq c_2$.
\end{lemma}
\proof From Zolotarev (1986), Th. $2.2.3.$, formula $(2.2.18)$, we get
\begin{equation} \label{A8.1}
p(x)=\frac{r|x|^{1/(r-1)}}{2 (r-1)} \int_0^1 u(\varphi)
\exp(-|x|^{r/(r-1)}u(\varphi)) d\varphi, \quad x \ne 0,
\end{equation}
where
\[
u(\varphi)=\left( \frac{\sin (\pi r \varphi/2)}{\cos (\pi
\varphi/2)}\right)^ {r/(1-r)} \frac{\cos (\pi
(r-1)\varphi/2)}{\cos(\pi \varphi/2)}.
\]
%(We have $\alpha=r$, $\beta=0$ in the notation of Zolotarev, then $\theta=0$ and
%$\theta^*=0$ in his notation).
Clearly, for $\varphi \in [1/2,1]$ we have
\begin{eqnarray*}
& & 1 \geq \cos (\pi (r-1)\varphi/2) \geq \cos (\pi (r-1)/4) >0\\
& & c_3\geq \sin(\pi r \varphi/2) \geq c_4 >0,
\end{eqnarray*}
where $c_3>0$ and $c_4>0$ are constants. Thus,
\[
c_6 \left( \cos (\pi \varphi/2)\right)^{1/(r-1)} \leq u(\varphi) \leq
c_5 \left( \cos (\pi \varphi/2)\right)^{1/(r-1)},
\]
$\varphi \in [1/2,1]$, $c_5>0$, $c_6>0$ are constants. Now, if $\vv \in [1/2,1]$
\[
c_7 (1 - \varphi) \leq \cos (\pi \varphi/2) \leq c_8 (1 - \varphi)
\]
for some $c_7>0$, $c_8>0$. Finally,
\[
c_{10} (1- \varphi)^{1/(r-1)} \leq u(\varphi) \leq c_9 (1- \varphi)^{1/(r-1)},
\forall \varphi \in [1/2,1]. \label{L1}
\]
Using $(\ref{A8.1})$ and the fact that $u(\varphi) \geq 0$ for $\varphi \in [0,1]$,
we get
\begin{eqnarray*}
p(x) &\geq & c |x|^{1/(r-1)} \int_{1/2}^1 (1 - \varphi)^{1/(r-1)}
\exp \left( -|x|^{r/(r-1)} c_9 (1- \varphi)^{1/(r-1)} \right) d \varphi\\
&=& c |x|^{1/(r-1)} \int_0^{1/2} \varphi^{1/(r-1)}
\exp \left( - c_9 (|x|^r \varphi)^{1/(r-1)} \right) d \varphi.
\end{eqnarray*}
Here and further on $c>0$ are constants, probably different on different occasions.

By change of variables, $u=(|x|^r \varphi)^{1/(r-1)}$, we get
\begin{eqnarray*}
p(x) & \geq & c |x|^{1/(r-1)} \int_0^{(|x|^r/2)^{1/(r-1)}} \frac
{u}{|x|^{r/(r-1)}} \exp(-c_9 u) \frac{u^{r-2}}{|x|^r} du\\
& = &  c |x|^{-1-r} \int_0^{(|x|^r/2)^{1/(r-1)}} u^{r-1} \exp(-c_9 u)du\\
& \geq & c |x|^{-1-r} \int_0^{(c_2^r/2)^{1/(r-1)}} u^{r-1} \exp(-c_9 u)du
\geq c_1 |x|^{-1-r},
\end{eqnarray*}
for $|x| \geq c_2 >0$. This also implies that $p(x)>0$, $\forall x\ne 0$, and
$$
p(0) = (2 \pi)^{-1} \int \exp (-|t|^r) dt \ne 0,
$$
hence $p$ is positive on $\mathbb{R}$.
\hfill $\Box$

%%%%%%%%%%%%%%%%%%%%%%%%%%%%%%%%%%%%%%%%
\begin{lemma} \label{nine}
Let $0 < r < s < \infty$ and let $h_* = h_*(n)$ be defined by
$(\ref{hstar})$, i.e.
\begin{eqnarray*}
  \frac{2 \alpha}{h_*^r}+\frac{2 \beta}{h_*^s} &=& \log n -(\log \log n)^2.
\end{eqnarray*}
Let $h_n$ satisfy
\begin{eqnarray*}
  b \log h_n + \frac{2 \alpha}{h_n^r}+\frac{2 \beta}{h_n^s} &=& \log n + C(1+o(1)),
 \quad n \to \infty,
\end{eqnarray*}
for some $b \in \r$ and $C \in \r$. Then, as $n\to \infty$, we have
\begin{equation}\label{AA}
h_*(n)
=(\log n/(2 \beta))^{-1/s}(1+o(1)),
\end{equation}
\begin{equation}\label{A}
    h_n^{a} \exp\left(-\frac{2 \alpha}{h_n^r} \right)
    = h_*^{a} \exp\left(-\frac{2 \alpha}{h_*^r} \right)(1+o(1)),
\end{equation}
\begin{equation}\label{B}
    \frac{h_*^{a}}{n} \exp\left(\frac{2 \beta}{h_*^s} \right)
    = o\left(\exp\left(-\frac{2 \alpha}{h_*^r} \right)\right),
\end{equation}
for any $a\in\r$, and
\begin{equation}\label{C}
    h_*^{s+2\gamma-1} \exp\left(\frac{2 \beta}{h_*^s} \right)
    \leq h_n^{s+2\gamma-1} \exp\left(\frac{2 \beta}{h_n^s} \right),
\end{equation}
for $n$ large enough.
\end{lemma}
\proof
Define $x_*=h_*^{-s}$, $x_n=h_n^{-s}$, and write, for $t>0$,
$$
F(t)\egal 2 \beta t + 2 \alpha t^{r/s}, \quad
F_1(t) \egal (-b/s) \log t +2\beta t +2\alpha t^{r/s}.
$$
Then
\begin{eqnarray}
  F(x_*) &=& \log n -(\log \log n)^2,
  \label{F1} \\
  F_1(x_n) &=&
  \log n + C(1+o(1)), \label{F2}
\end{eqnarray}
for a constant $C\in \r$. We first prove that $x_n$ satisfies
\begin{equation}\label{F3}
    F(x_n) =\log n+C_1 \log \log n(1+o(1))
    +C_2(1+o(1))
\end{equation}
for some constants $C_1,C_2\in \r$. In fact,
\begin{eqnarray*}
F_1^\prime (x_n) &=& \frac{1}{x_n}\left(-\frac{b}{s} \right) +2\beta
+ \frac{2\alpha r}{s}x_n^{r/s-1}>0,
\end{eqnarray*}
for $x_n$ large enough, thus $F_1(t)$ is strictly monotone increasing for large $t$,
and a solution $x_n$ of $(\ref{F2})$ exists for large $n$ (and is unique).
Next, clearly,
$$
\frac{F_1(t)}{2 \beta t} \to 1,\quad  t \to \infty,
$$
and therefore $\log n/(2 \beta x_n) \to 1$, as $n \to \infty$. Similarly,
$\log n/(2 \beta x_*) \to 1$, as $n \to \infty$, which yields (\ref{AA}). Thus
$(-b/s) \log x_n = (-b/s) \log \log n (1+o(1))$, as $n \to \infty$, and
write $F(x_n) = F_1(x_n) +(b/s) \log x_n$ to get $(\ref{F3})$ in view of
$(\ref{F2})$.
We have
$$
x_n=F^{-1}(\log n +a_n),\quad x_*=F^{-1}(\log n-b_n)
$$
where $a_n=C_1 \log\log n(1+o(1))+C_2(1+o(1))=O(\log \log n)$,
$b_n=(\log \log n)^2$ and $F^{-1}(\cdot)$ is the inverse of $F(\cdot)$.
Hence, for some $0< \tau <1$ and for $n$ large enough,
\begin{eqnarray}
x_n &=& F^{-1}(\log n+a_n)= x_* +(F^{-1}(\log n-b_n))^\prime(a_n+b_n)\nonumber \\
&&+\frac{1}{2} \left(F^{-1}(\log n-b_n(1-\tau)+ \tau a_n)\right)''(a_n+b_n)^2.
\label{ff1}
\end{eqnarray}
The first and the second derivatives of $F^{-1}$ are given by
\begin{eqnarray*}
(F^{-1}(y))^\prime &=& \frac{1}{F^\prime(F^{-1}(y))}
= \frac{1}{2 \beta + (2 \alpha r/s)(F^{-1}(y))^{r/s-1}},\\
(F^{-1}(y))^{\prime \prime} &=& \frac{-(2 \alpha r/s)(r/s - 1)(F^{-1}(y))^{r/s-2}}
{(2 \beta +(2 \alpha r/s)(F^{-1}(y))^{r/s - 1})^3}.
\end{eqnarray*}
Hence
\begin{equation}
\label{ff2}
(F^{-1}(\log n-b_n))^\prime = \frac{1}{2 \beta +(2 \alpha r/s) x_*^{r/s-1}}
= \frac{1}{2\beta}+o(1), \quad n \to \infty.
\end{equation}
Next, it easy to show
that there exists $\bar{y}>0$ such that
\begin{equation}
\label{ff3}
y/(4 \beta) \le F^{-1}(y) \leq y/(2 \beta)
\end{equation}
for $y \geq \bar{y}$.
Considering $n$ large enough so that $y_n\egal \log n - b_n(1-\tau)+a_n \tau \geq \bar{y}$
and using the above expression for $(F^{-1}(y))^{\prime \prime}$ and (\ref{ff3}) we get
\begin{eqnarray*}
(F^{-1}(y_n))^{\prime \prime}&=&
(F^{-1}(\log n - b_n(1-\tau)+ \tau a_n))^{\prime \prime}= O((\log n)^{r/s-2}),
\quad n \to \infty.
\end{eqnarray*}
This and (\ref{ff1}), (\ref{ff2}) imply
\begin{equation}
  x_*-x_n = -\frac{1}{2 \beta}(1+o(1)) (a_n+b_n)+O\left(\frac{(a_n+b_n)^2}
  {(\log n)^{2-r/s}} \right)
  %\nonumber\\
 % &=& -\frac{b_n}{2 \beta} (1+o(1)) + O\left(\frac{(\log \log n)^4}
 %{(\log n)^{2-r/s}} \right)
  %\nonumber\\
  %&=&
  = -\frac{b_n}{2 \beta} (1+o(1)).
  \label{posl}
\end{equation}
Using this representation we obtain
\begin{eqnarray*}
  && \exp\left(-\frac{2 \alpha}{h_n^r}+\frac{2 \alpha}{h_*^r} \right)
  = \exp(-2\alpha (x_n^{r/s}-x_*^{r/s})) \\
  %&=& \exp\left( -2\alpha [x_*+b_n(2\beta)^{-1}(1+o(1))]^{r/s}-x_*^{r/s}) \right) \\
  &=& \exp\left( -2\alpha x_*^{r/s}([1+b_n(2\beta x_*)^{-1}(1+o(1))]^{r/s}-1)
  \right) \\
  &=& \exp\left(O(b_n x_*^{r/s-1}) \right) = 1+o(1),
\end{eqnarray*}
since $b_n=(\log \log n)^2$, $x_* = (2 \beta)^{-1} \log n \ (1+o(1))$ and $r<s$.
This and the fact that $(h_n / h_*)^{a} = (x_*/x_n)^{a/s}= 1+o(1)$
imply $(\ref{A})$. Next, (\ref{B}) follows directly from the definition of
$h_*$ and from (\ref{AA}). To prove $(\ref{C})$, note that, in view of (\ref{posl}),
\begin{eqnarray*}
  && \frac{h_*^{s+2\gamma-1}}{h_n^{s+2\gamma-1}} \exp\left(\frac{2 \beta}{h_*^s} -
  \frac{2 \beta}{h_n^s} \right)=(1+o(1))\exp(2\beta(x_*-x_n)) \\
  &=& (1+o(1)) \exp(-b_n[1+o(1)]) \leq 1
\end{eqnarray*}
for $n$ large enough.
\hfill $\Box$

\begin{lemma} \label{two}
Let $0 < r < s < \infty$ and let $h_+ = h_+(n)$ be the solution of
$(\ref{LB6})$.
Then $h_+(n) =(\log n/(2 \beta))^{-1/s}(1+o(1))$,
\begin{equation}\label{2L.1}
    h_+^{a} \exp\left(-\frac{2 \alpha}{h_+^r} \right)
    = h_*^{a} \exp\left(-\frac{2 \alpha}{h_*^r} \right)(1+o(1)), \text{ as }n\to \infty
\end{equation}
and
\begin{equation}\label{2L.2}
    (\log n)^b n \exp \left(- \frac{2 \alpha}{h_+^r}-\frac{2 \beta}{h_+^s}\right) =o(1),
\end{equation}
as $n \to \infty$, for any $a \in \r$, $b \in \r$.
\end{lemma}

\noindent{\bf Proof} is analogous to that of Lemma~\ref{nine}. \hfill $\Box$

\bigskip

\noindent{\bf Acknowledgement.}  The results of  this paper were  presented at
the Conference ``Rencontres de statistiques math\'ematiques'', CIRM Luminy, 2001.
Later, Fabienne Comte and Marie-Luce Taupin suggested a different estimator
for the same problem refraining from studying the optimality of rates issue
(Comte and Taupin (2003)). We would like to thank them for discussion of the results.

\footnotesize{

}
\begin{flushright}
\footnotesize
$^1${\sc Laboratoire de Probabilit\'es et Mod\`eles }\\
{\sc Al\'eatoires (UMR CNRS 7599),}\\
{\sc Universit\'e Paris VI}\\
{\sc 4, pl.Jussieu, Bo\^{\i}te courrier 188,}\\
{\sc 75252 Paris, France}\\
{\sc e-mail}: tsybakov@ccr.jussieu.fr
\bigskip\\
$^2${\sc Modal'X, Universit\'e Paris X}\\
{\sc 200, avenue de la R\'epublique }\\
{\sc 92001 Nanterre Cedex, France}\\
{\sc e-mail}: butucea@ccr.jussieu.fr

\end{flushright}

\end{document}